\documentclass[11pt]{article} %

\usepackage[preprint]{tmlr}

\usepackage{mathpazo}

\usepackage{microtype}

\usepackage{amsmath,amsfonts,bm}

\def\eqref#1{equation~\ref{#1}}

\def\1{\bm{1}}

\def\mM{{\bm{M}}}

\def\mR{{\bm{R}}}

\DeclareMathAlphabet{\mathsfit}{\encodingdefault}{\sfdefault}{m}{sl}
\SetMathAlphabet{\mathsfit}{bold}{\encodingdefault}{\sfdefault}{bx}{n}

\usepackage{amsmath}
\usepackage{amsfonts} 
\usepackage{mathtools}
\usepackage{amsthm} %
\usepackage{amssymb}
\usepackage{mathabx}

\usepackage{graphicx}
\usepackage{subcaption}

\usepackage{graphicx}
\usepackage{subcaption}
\usepackage{amsthm} %
\usepackage[e]{esvect}
\usepackage{enumerate}
\usepackage{enumitem}

\usepackage{hyperref}
\usepackage{url}

\usepackage{xcolor}
\usepackage{mdframed}

\definecolor{mygray}{gray}{0.92}
\definecolor{mycolor2}{RGB}{255, 250, 176}

\newmdenv[backgroundcolor=mygray, linewidth=0pt]{sectionpurpose}

\usepackage{xspace}
\newcommand*{\eg}{e.g.\@\xspace}
\newcommand*{\ie}{i.e.\@\xspace}

\newcommand*{\cf}{cf.\@\xspace}
\newcommand*{\versus}{vs.\@\xspace}

\newtheorem{theorem}{Theorem}[section]       %
\newtheorem{proposition}[theorem]{Proposition}
\newtheorem{definition}[theorem]{Definition}
\newtheorem{corollary}[theorem]{Corollary}
\newtheorem{lemma}[theorem]{Lemma}
\newtheorem{remark}[theorem]{Remark}
\newtheorem{example}[theorem]{Example}
\newtheorem{conjecture}[theorem]{Conjecture}

\usepackage{xspace}

\newcommand{\RomanNumeralCaps}[1]
    {\MakeUppercase{\romannumeral #1}}

\renewcommand{\mM}{\mathcal{M}}

\renewcommand{\mR}{\mathbb{R}}
\newcommand{\Pup}{\overline{P}}
\newcommand{\Plow}{\underline{P}}
\newcommand{\Rup}{\overline{R}}
\newcommand{\Rlow}{\underline{R}}

\newcommand{\CP}{\operatorname{CP}}
\newcommand{\limsupn}{\limsup_{n \rightarrow \infty}}
\newcommand{\limn}{\lim_{n\rightarrow \infty}}

\newcommand{\linfty}{L^\infty}
\newcommand{\linftyd}{(\linfty)^*}

\newcommand{\cobar}{\bar{\operatorname{co}}}
\newcommand{\co}{\operatorname{co}}
\newcommand{\sumin}{\sum_{i=1}^n}

\newcommand{\ccm}{\oplus}
\newcommand{\iid}{\textit{i.i.d.}\@\xspace}
\newcommand{\etc}{etc.\@\xspace}

\newcommand{\PlowWF}{\Plow^\infty_{\mM,\operatorname{WF}}}

\newcommand{\NSLP}{NSLP}
\newcommand{\SLI}{SLI}
\newcommand{\infprodsigma}{\mathcal{B}(\Omega^\infty)}

\newcommand{\avgest}{\hat{\underline{\operatorname{avg}}}_{\ell,n}}
\newcommand{\avgesta}{\hat{\underline{\operatorname{avg}}}_{\chi_{A},n}}

\makeatletter
\newcommand{\shorteq}{%
  \settowidth{\@tempdima}{-}%
  \resizebox{\@tempdima}{\height}{=}%
}
\makeatother
\newcommand{\DtoN}{\text{CP} \shorteq \hspace{0.3mm}\mathcal{N}}

\newcommand{\MtoNset}{\mM_{\shortrightarrow \mathcal{N}}^\infty}
\newcommand{\CM}{\mathcal{C}(\mM)}
\newcommand{\kfunc}{\tilde{k}}

\usepackage{stmaryrd}

\title{Data Models With Two Manifestations of Imprecision}

\author{\name Christian Fröhlich \email christian.froehlich@uni-tuebingen.de \\
      \addr Department of Computer Science\\
      University of Tübingen\\
      and Tübingen AI Center
      \AND
      \name Robert C. Williamson \email bob.williamson@uni-tuebingen.de \\
      \addr Department of Computer Science\\
      University of Tübingen\\ 
      and Tübingen AI Center
      }

\def\month{MM}  %
\def\year{YYYY} %
\def\openreview{\url{https://openreview.net/forum?id=XXXX}} %

\begin{document}

\maketitle

\begin{abstract}
Motivated by recently emerging problems in machine learning and statistics, we propose data models which relax the familiar \iid assumption. 
In essence, we seek to understand what it means for data to come from a set of probability measures. We show that our frequentist data models, parameterized by such sets,  manifest two aspects of imprecision. 
We characterize the intricate interplay of these manifestations, aggregate (ir)regularity and local (ir)regularity, where a much richer set of behaviours compared to an \iid model is possible. In doing so we shed new light on the relationship between non-stationary, locally precise and stationary, locally imprecise data models. We discuss possible applications of these data models in machine learning and how the set of probabilities can be estimated. For the estimation of aggregate irregularity, we provide a negative result but argue that it does not warrant pessimism. Understanding these frequentist aspects of imprecise probabilities paves the way for deriving generalization of proper scoring rules and calibration to the imprecise case, which can then contribute to tackling practical problems.
\end{abstract}

\section{Introduction}
In machine learning and statistics, the assumption that data is drawn \iid is ubiquitous, so firmly entrenched to the extent that it is seldom questioned. Recently, however, various challenges to the \iid assumption have surfaced and received increasing attention, among them the problem settings of \textit{data corruptions} (including \textit{dataset shift}) \citep{iacovissi2023general}, \textit{fairness desiderata}, \textit{federated learning}, \textit{multi-source adaptation} and others. These settings call for a relaxation of the \iid assumption to the weaker assumption that the data comes from a \textit{set of probabilities} instead of a single one.

Under the banner of \textit{imprecise probabilities} (IP; \citet{walley1991statistical}) scholars have gathered to systematically study more flexible models, essentially equivalent to sets of probability measures. With relatively few exceptions, work on IP has been rooted in a subjectivist approach, taking a set of probability measures to represent a state of \textit{belief}, which justifies the name \textit{credal set} for such a set of probability measures. As a connection to the frequentist side, an often-used interpretation for imprecise probabilities is that data is drawn \iid, but since the true precise probability is not known to the agent, the imprecision supposedly models the \textit{epistemic uncertainty} arising from the lack of knowledge (see \eg \citep{hullermeier2022quantification}). This is in line with some works in the paradigm of \textit{distributional robustness}, which assume a true but unknown precise probability and aim to construct an ambiguity set of probability measures to guard against the ``estimation'' error \citep{kuhn2019wasserstein}. 
In this view, knowledge of the true probability eliminates imprecision. In contrast, the problems listed above which motivate us (data corruptions \etc) are not appropriately modelled as \iid, even when putting epistemic concerns aside. Few works in the IP literature \citep{walley1982towards} have taken an interpretation of \textit{aleatoric indeterminacy} or what is also called an \textit{ontic interpretation} \citep{augustin2022statistics}, which grounds the imprecision not in a lack of knowledge of the data-generating process, but in frequentist aspects of the data itself: the data is drawn from a set of probabilities --- or rather, it is modelled in this way.

Taking this perspective, in this paper we relax the \iid assumption.
First, we consider a data model which still assumes precise, probabilistic independence, but where we allow a set of probability measures instead of a single one. We call this the \textit{non-stationary, locally precise} data model.\footnote{Of course we are not the first ones to use such a basic model, but we do study the aspects of aggregate (ir)regularity and local (ir)regularity of these models in detail.} It is specified by an equation of the form (stated rigorously in Section~\ref{sec:nonstationarymodels})
\[
\lambda\left\{ \omega^\infty = (\omega_1, \omega_2, ..)\right\} = \prod_{i=1}^\infty p_i(\{\omega_i\}),
\]
where $\omega^\infty$ is an infinite sequence of elementary events, and the $p_i$ are probability measures. This models is locally precise in the sense that for each $i \in \mathbb{N}$, a precise probability $p_i$ governs the data generation. The imprecision is thus, in contrast, a global feature of the model. We are then interested in studying the implications of such a \textit{data model with underlying (global) imprecision}.
In particular, we focus on two aspects that the familiar \iid model exhibits, \textit{aggregate regularity} and \textit{local irregularity} (randomness), and our main theorem fully characterizes how these aspects are transformed and interact in such more general data models (Theorem~\ref{theorem:maintheoremffivanenko}). In summary:

\begin{itemize}
    \item Under an \iid model, relative frequencies of events converge \textit{almost surely} to a limit. In our data models with underlying imprecision, this might break down. Here, the natural generalization is to consider the \textit{set of cluster points} (accumulation points) of relative frequencies.
    \item Under an \iid model, outcome sequences are \textit{almost surely} ``perfectly random'' (to be explained later). In our data models with underlying imprecision, we may find hidden heterogeneity, that is, subsequences with relative frequencies that even in the limit do not coincide with the aggregate regularity. For example, in a fairness context, this would correspond to subpopulations for which the probability is distinct from the population probability.
\end{itemize}

Previous frequentist work in the literature can be broadly distinguished by which of these two \textit{manifestations} of imprecision they focus on. \citet{walley1982towards} and \citet{ivanenkobook} have focused on the first, whereas a later line of work by Fine and colleagues (\eg \citet{fierens2009frequentist}) studies the second manifestation. These works can be, in our opinion, challenging to access. In our paper, we aim to provide an easily accessible account of frequentist imprecision, which also for the first time explicates the relation between these two manifestations and thereby unifies previous works.

In the process of establishing relations to a line of work initiated by \citet{walley1982towards} (Section~\ref{sec:comparisontofine}), we go beyond non-stationarity and local precision. By using certain non-stationary, locally precise data models as basic building blocks, we arrive at more general data models, which are stationary and where the imprecision is now already a local feature, and sometimes a global feature. This puts us in a position to also compare our approach (Section~\ref{sec:llncomparisonsection}) to a literature on generalized laws of large numbers for imprecise probabilities and coherent risk measures, a strain of literature which is actively being developed (\eg \citep{hu2016general,peng2019nonlinear,zhang2024conditional}). Here, the distinction of \textit{risk} and \textit{ambiguity} is a prominent motivation. This comparison also illustrates how data models relate to subjectivist approaches more generally.

Finally, to the possible implications of this work. We hope that our paper can provide a starting point for addressing practical problems outlined above (and others) in a more principled manner by providing a conceptual foundation. In the classical picture, the \iid data model is paired with \textit{proper scoring rules} and \textit{calibration} to obtain a general framework for obtaining and evaluating probabilistic forecasts. We believe that introducing generalized proper scoring rules and calibration, adapted to the imprecise case, is required, and can then lead to the development  of practically useful methods. Our position is that these concepts need to be considered as relative to a data model, and in this way our work makes first steps in this direction. 
Developing these concepts would also serve a related goal: 
in the philosophical literature on rational decision making there is a heated debate about whether imprecise belief state are ever warranted (\eg \citep{elga2010subjective,bradley2014should,schoenfield2017accuracy})); for example, when might a rational decision maker justifiably believe that the probability of rain tomorrow is between $20\%$ and $40\%$? Our approach makes it plausible that one part of the controversy could be resolved: at least under a data model with underlying imprecision, imprecise belief states are warranted. Still, in our opinion, a definite answers will require conceptually developing proper scoring rules and calibration for such imprecise data models. Also, this reasoning remains silent on the case for imprecise belief states when the data model is precise and the goal is to capture epistemic uncertainty with imprecision, as in the paradigm of distributional robustness.
While developing proper scoring rules and calibration for imprecise data models is out of scope for the present paper,
we show how the two manifestations of imprecision can (and \textit{cannot}) in principle be estimated.

The paper is organized as follows. We begin by a brief, general treatment of \textit{data models} to set the stage for going beyond the \iid model (Section~\ref{sec:datamodels}). We highlight two aspects of the \iid model: aggregate regularity and local irregularity. We discuss a further crucial ingredient in any data model, \textit{typicality}. We introduce the non-stationary, locally precise data model based on a set of probabilities (Section~\ref{sec:nonstationarymodels}). Our main theorem characterizes possible manifestations of imprecision in such models as aggregate (ir)regularity and local (ir)regularity (Theorem~\ref{theorem:maintheoremffivanenko}). We introduce stationary, locally imprecise data models in Section~\ref{sec:comparisontofine}, compare them to the highly related work of \citet{walley1982towards} and subsequent papers, and to generalized laws of large numbers (Section~\ref{sec:llncomparisonsection}).
We discuss possible applications for our data models (Section~\ref{sec:applications}) and make some general remarks about estimation. Here, we provide a negative result concerning an estimator proposed by \citet{walley1982towards}, but argue that it is not as troubling as it may appear at first sight. 
We conclude by motivating the need for imprecise scoring rules and imprecise calibration, relativized to such data models.

\textit{Who should read this paper?} We believe this paper might be of relevance to machine learning and statistics scholars interested in establishing conceptual foundations for the above practical (and similar) problems, which escape being modelled as \iid, as well as scholars generally interested in frequentist aspects of imprecise probabilities and coherent risk measures. In particular, those who have found previous works on such frequentist aspects hard to access (\eg \citep{walley1982towards}) 
 may find our exposition helpful.

\textit{General notation}. For a set $A$, we denote its powerset as $2^A$. We write $\cobar(A)$ for the closed convex hull in the appropriate topology. By $\chi_A$ we denote the indicator function of the set $A$, defined as $\chi_A(\omega) \coloneqq 1$ if $\omega \in A$, $0$ otherwise. We write the set of $A$-valued sequences as $A^\infty$. An element of $A^\infty$ is some infinite sequence $a^\infty=(a_1,a_2,..) \in A^\infty$. For a finite $n \in \mathbb{N}$, we define $A^n \coloneqq \{a=(a_1,..,a_n) : a_i \in A\}$.
If $a=(a_1,..,a_u) \in A^u$ and $b=(b_1,..,b_v) \in A^v$, if $u \leq v$, we write $a \subseteq b$ to mean that $a$ is a prefix of $b$, and likewise if $b$ is an infinite sequence. We use the terms \textit{probability} and \textit{probability measure} synonymously, but sometimes choose the latter to stress the measure-theoretic character.

\section{Data Models}
\label{sec:datamodels}
Assume a finite \textit{possibility set} $\Omega = \{\omega^1, .., \omega^k\}$, where an elementary event $\omega^i \in \Omega$ captures all relevant aspects of the world given a decision-making context. Throughout the paper, we will assume finite $\Omega$ and $1<|\Omega|=k$. 
In line with the IP literature, we call bounded functions $X: \Omega \rightarrow \mR$ \textit{gambles} and collect them in the set $\linfty$. For example, $\Omega = \mathcal{X} \times \mathcal{Y}$ might be a joint space containing features and labels, which are then extracted by  $X,Y \in \linfty$; since $\mR$ and $\mR^d$ have the same cardinality, this allows handling vector-valued features. Given a loss function $\ell : \mathcal{A} \times \Omega \rightarrow \mR$, we are typically interested in the gamble $\ell(a,\cdot)$ for a chosen action $a \in \mathcal{A}$.\footnote{The IP literature usually works with a gain orientation, where positive values correspond to gain (utility), whereas we interpret positive values as losses (disutility). No difference results beyond simple sign flips.}

By \textit{data} (or \textit{data sequences}) we understand a finite or infinite sequence of elementary events $\omega^\infty = (\omega_1, \omega_2, ..)$. In the following, we restrict ourselves to infinite sequences and denote the set of such sequences as $\Omega^\infty$.
Data is that which is actually fixed. In contrast, a \textit{data model} is a description of a data-generating process, for specifying what kind of data we expect to observe.
For a data model we require a \textit{generator}, which is a function $G : \mathcal{Z} \to \Omega^\infty$ which takes in a \textit{seed} $z \in \mathcal{Z}$ and returns data. Such a seed can be understood as an abstract source of \textit{process randomness} (which we contrast below with \textit{outcome randomness}). 
A standard choice is to use $\mathcal{Z}=[0,1)$, since a real number from the unit interval is via its binary expansion identified with an infinite binary sequence \citep{williams1991probability}. To make this more concrete, consider the familiar \iid model, though perhaps expressed here in an unfamiliar way. Denote the Lebesgue measure as $\lambda$. 

\begin{definition}[The \iid model]
    Let $\mathcal{Z} \coloneqq [0,1)$. Given a probability measure $p$ on $\Omega$, it is possible to consistently define $W_i : [0,1) \rightarrow \Omega$, $i \in \mathbb{N}$ so that the following holds:
    \[
    \label{eq:iidmodeleq}
\lambda\left\{z \in \mathcal{Z} : W_1(z) = \omega_1, W_2(z)=\omega_2, .., W_n(z)=\omega_n \right\} \coloneqq \prod_{i=1}^n p(\{\omega_i\}), \quad \forall n \in \mathbb{N},
\]
Define $G : [0,1) \to \Omega^\infty$ as $G(z) \coloneqq (W_1(z),W_2(z),..)$.
\end{definition}
Note that the $W_i$ are the quantities that are being defined here, not $p$ or $\lambda$. For a proof that the $W_i$ can be defined so that \eqref{eq:iidmodeleq} holds, see \citep[Theorem 3.16]{spreij2012measure}. The \iid model forms the departure point for our more general data models and it is thus worth highlighting some of its features. We do so in a way that is grounded in Richard von Mises \citeyearpar{mises1919grundlagen} foundational account of probability.

\textbf{Aggregate regularity}, also called \textit{statistical stability} \citep{gorban2017statistical}. 
Under an \iid model it holds that
\begin{equation}
    \label{eq:aggreg}
    \lambda\left\{z \in \mathcal{Z} : G(z) = (\omega_1, \omega_2, ..) : \forall A \subseteq \Omega: \lim_{n \rightarrow \infty} \frac{1}{n} \sum_{i=1}^n \chi_A(\omega_i) = p(A)\right\} = 1.
\end{equation}        
 We call $r_n(A) \coloneqq \frac{1}{n} \sum_{i=1}^n \chi_A(\omega_i)$ the relative frequency of the event $A$. In words, ``almost surely, relative frequencies converge to their \textit{probability}''. While the above statement, the \textit{strong law of large numbers}, is of formal, mathematical nature, there exists an associated empirical sibling, what \citet{gorban2017statistical} calls the \textit{hypothesis of statistical stability}. This refers to the statistician's in principle unfalsifiable belief that the data at hand is the prefix of a statistically stable sequence, that is, gathering more and more data will lead to stabilization of relative frequencies and gamble averages. Note that this belief is not logically implied by the theoretical law of large numbers, which is a formal property of the \iid data model. Already \citet[p.\@\xspace 175]{de2017theory} has noted
 \begin{quote}
     ``In any case, the force of the ‘stability of frequencies’ as a probabilistic or
statistical principle is completely illusory, and without solid foundation.''
 \end{quote}
 
 In fact, \citet{gorban2017statistical} has demonstrated empirically that one does not always observe stable relative frequencies even for very long observation intervals. Recent studies in climate science also observe this \citep{lovejoy2015voyage,franzke2020structure}. Consequently, it seems prudent to explore more flexible data models. We also refer the reader to \citep{frohlich2024strictly} for a discussion of statistical stability.

\textbf{Local irregularity}, also called \textit{randomness}. Under an \iid model, the regularity emerges only at the aggregate level, while the individual level remains fully unpredictable (random), relative to the probability $p$.\footnote{Intuitively, the best forecast when being judged with a proper scoring rule that a forecaster could announce is the true probability.} One way of mathematizing this intuition is with von Mises \citeyearpar{mises1919grundlagen} concept of \textit{selection rules}. For our purposes, a selection rules is a function $S : \mathbb{N} \to \{0,1\}$, which either selects or does not select at a given index $i \in \mathbb{N}$. Such a selection rule effectively extracts a subsequence. Then, for a countable set of selection rules $\mathcal{S}$, under the \iid model $(\forall A \subseteq \Omega)$:
\begin{equation}
    \label{eq:localirreg}
     \lambda\left\{z \in \mathcal{Z} : G(z) = (\omega_1, \omega_2, ..) : \forall S \in \mathcal{S}: \lim_{n \rightarrow \infty} \frac{\sum_{i=1}^n \chi_A(\omega_i) S(i)}{\sum_{i=1}^n S(i)} = p(A)\right\} = 1.
\end{equation}      
That limiting relative frequencies along a subsequence coincide with the probability as specified by $p$ is a form of randomness desideratum. Intuitively, a skeptic who uses a selection rule $S$ to decide whether or not to place a bet of the form $\omega_i - p(A)$ on the next outcome cannot get infinitely rich, when accumulating capital as $i \to \infty$. In light of the above assertion, we can expect ``perfect randomness'' under an \iid model. By contrast, assume that \eqref{eq:localirreg} would \textit{not} hold, meaning some $S \in \mathcal{S}$ reveals what we might call \textit{local regularity}, \textit{non-randomness} or \textit{hidden heterogeneity}. In a fairness context, for example, this would mean that the probability associated to a subpopulation does not coincide with the probability associated to the whole population (see \citep[p.\@\xspace 175]{de2017theory}). In the context of machine learning, optimizing only an aggregate criterion (\eg average loss) then implies risking unfairness (non-robustness) \citep{williamson2019fairness}. The term \textit{hidden heterogeneity} highlights that since all information about a single datum is captured by the elementary event $\omega$ (recall that the feature map is $X: \Omega \to \mR$), a predictor $f(X)$ cannot, if there are distinct probabilities on $\Omega$, be ``aware'' of this heterogeneity. This also illustrates that whether there is hidden heterogeneity depends on the modelling, the choice of features. If for example, in a fairness context, the sensitive feature is added to $\mathcal{X}$, the heterogeneity is ``revealed'' and a predictor can make use of it. In this way, hidden heterogeneity refers to a non-randomness that is not captured by the probabilistic randomness of the features.

We remark that, in this paper, we brush aside a subtle distinction: in von Mises framework, a selection rule is a function that can depend on all previous outcomes (meaning a skeptic could take previous outcomes into account when deciding whether to place a bet or not), whereas in our setting it may only depend on the index $i \in \mathbb{N}$. For a fixed sequence, the formulations coincide. Mathematically, our setting is significantly simpler, however.

While the two features above are specific to the \iid model (we will later ``break'' them), another ingredient of any data model in our view is a notion of \textbf{typicality}. Some data sequences possess specific properties. For instance, we might be interested in sequences with converging relative frequencies. Such a property can then be identified with a subset of seeds:
\[
\mathcal{Z}_A  \coloneqq \{z \in \mathcal{Z}: G(z) \text{ has property } A\} \subseteq \mathcal{Z}.
\]
Both statements \eqref{eq:aggreg} and \eqref{eq:localirreg} hold with the qualification that the Lebesgue measure of the set of seeds which generate sequences with the property of interest is $1$; not that it holds \textit{for all} sequences.
We think the essence of such an assertion lies in asserting \textit{typicality}: it means that \textit{typically} we will observe sequences with converging relatives frequencies; other sequences are so rare that they are considered negligible. An advantage of framing this in terms of typicality rather than probability, is that it circumvents the following problematic move by frequentists, discussed by \citet{lacazefrequentism}:
\begin{quote}
    ``Importantly, `almost sure convergence' is also given a frequentist interpretation. Almost
sure convergence is taken to provide a justification for assuming that the relative frequency
of an attribute \textit{would} converge to the probability in actual experiments \textit{were} the experiment
to be repeated indefinitely'' [emphasis in original].
\end{quote}
That is, there are actually two measures at play, but one expresses typicality, whereas the other expresses probability; they are conjoined, indeed, but only \textit{typically}.
For our purposes, we do not require interpreting the typicality statement as a frequentist statement itself, we simply accept the more parsimonious notion that typicality is somehow defined and can be used instrumentally (that it has some justification). 
When we speak of \textit{probability}, we mean a measure on $\Omega$ which is appropriately normalized. In contrast,
the ``randomness'' in the data generation instead is tied to typicality and should be conceptually distinguished from what we have above called local irregularity, a property of a fixed data sequence; one way of phrasing this distinction is in terms of \textit{process} randomness (concerning typicality) \versus  \textit{outcome} (or \textit{product}) randomness (concerning local irregularity) \citep{sep-chance-randomness}. Thus, in our setup of the \iid model, the Lebesgue measure expresses typicality, we may call it \textit{probabilistic typicality} \citep{galvan2006bohmian} since the Lebesgue measure on $[0,1)$ formally is a probability measure in the sense of measure theory, but the semantics are distinct from that of a probability $p$ on $\Omega$. 
Conceptually distinguishing typicality from probability is valuable because it opens up room for \textit{nonprobabilistic typicality} notions which have been axiomatically introduced and studied by \citet{galvan2006bohmian}. In our case, we will use a \textit{coherent lower probability} as the typicality notion for our data models.

\citet{galvan2006bohmian} is, to our knowledge, the first to initiate an axiomatic study of typicality beyond probabilistic typicality. To this end, the author axiomatically defines the notion of a \textit{typicality distance} between sets of $\mathcal{Z}$ (in our context), giving rise also to a \textit{relative} and \textit{absolute typicality measure}. We refer the reader to \citet{galvan2006bohmian} for details,
and here simply remark that in the probabilistic case, the author suggests as an absolute typicality measure (adapted to our context):
\[
T_a(Z) \coloneqq \lambda(\Omega \setminus Z), \quad \forall Z \in \mathcal{B}([0,1)),
\]
where $\mathcal{B}([0,1))$ is the Borel $\sigma$-algebra on $[0,1)$. \citet{galvan2006bohmian} defines a \textit{typical set} as a set satisfying $T_a(Z) \approx 0$; to make this formally rigorous, we could demand that $T_a(Z)=0$, which suffices when we only consider the asymptotic case. Thus, under an \iid model, the set of sequences for which relative frequencies converge to $p$ forms a typical set. Note that different sets can be typical: for instance, the set of sequences for which relative frequencies converge to $p$ is a strict superset of the set of sequences for which relative frequencies converge to $p$ \textit{and} relative frequencies along a countable set of subsequences also converge to the same limit. Both are typical sets under an \iid model.\footnote{In Galvan's \citeyearpar{galvan2006bohmian} terminology, this idea can be expressed as relative typicality.} Which typical set we need to consider depends on what and how we try to estimate.

The notion of a typical set is crucial for estimation. An \iid model which is parameterized by a true precise probability $p$ \textit{typically} generates sequences where limiting relative frequencies coincide with $p$. This means that relative frequencies provide a means to estimate the data model's parameterization under the assumption that our data sequence belongs to a typical set. In the precise case, minimizing a proper scoring rule guarantees that we recover the true probability in the limit for typical sequences (see \citep{gneiting2007strictly}). We would like to extend this to the imprecise case, so that we find estimation procedures, which for typical sequences recover the underlying set of measures.

This concludes setting the stage for moving to sets of probabilities, where the associated typicality notion will also be imprecise. 

\subsection{Background on Imprecise Probabilities}
We keep the background on imprecise probabilities at a bare minimum and the refer the reader to \citep{walley1991statistical} for a comprehensive account, and to \citep{augustin2014introduction} for a more accessible introduction. For our purposes, an imprecise probability is essentially given as a set of probability measures. Assume a possibility set $\Omega$ with a set system $\mathcal{A}\subseteq 2^\Omega$, so that at least $\emptyset,\Omega \in \mathcal{A}$; often, $\mathcal{A}=2^\Omega$. In this introduction, to reduce the burden of technicalities, we assume $|\Omega|<\infty$. In Section~\ref{sec:stationarymodels}, we will need to consider imprecise probabilities on infinite possibility sets (namely a set of infinite sequences); we caution the reader that the statement below then hold only when replacing ``probability measures'' with finitely additive expectation functionals.\footnote{To be precise, with linear previsions $E \in \linftyd$, which are elements of the dual space of $\linfty$ endowed with the weak* topology, which satisfy $E(\chi_\Omega)=1$ and $E(X)\geq 0$ when $X \geq 0$. See \citep{walley1991statistical}. If $p$ is a countably additive probability measure, $\mathbb{E}_p[\cdot]$ is a finitely additive expectation, since countable additive implies finite additivity. In fact, all models considered in Section~\ref{sec:stationarymodels} arise as the envelopes of countably additive probability measures.}
\begin{definition}
\label{def:coherentlowerprob}
    We say that a set function $\Plow : \mathcal{A} \to [0,1]$ is a \textit{coherent lower probability} if
\[
\Plow(A) = \inf\{\mu(A) : \mu \in \mM\}, \quad A \subseteq \Omega,
\]
for some set $\mM$ of probability measures on $\Omega$.
\end{definition} 
A representation of this form is called \textit{envelope representation}.\footnote{Here, $\Plow$ is the \textit{envelope} of $\mM$. Other authors call the set $\mM$ itself the \textit{envelope}, \eg \citet{rockafellar2015measures}.}
The conjugate \textit{coherent upper probability} is given as $\Pup(A) \coloneqq 1-\Plow(A^c) = \sup\{\mu(A) : \mu \in \mM\}$. Clearly, for a conjugate pair of coherent lower and upper probability, $\Plow(A) \leq \Pup(A)$ $\forall A \subseteq \Omega$, in short $\Plow \leq \Pup$. Such a pair satisfies the following properties \citep[Section 2.7.4]{walley1991statistical} ($A,B \in \mathcal{A}$):
 \begin{enumerate}[nolistsep,start=1,label=\textbf{P\arabic*.}, ref=UP\arabic*]
  \item \label{item:p1} $\Plow(\emptyset)=\Pup(\emptyset)=0$, $\Plow(\Omega)=\Pup(\Omega)=1$.
  \item \label{item:p2} If $A \subseteq B$, then $\Plow(A)\leq \Plow(B)$ and $\Pup(A)\leq\Pup(B)$.
  \item \label{item:p3} $\Pup(A \cup B) \leq \Pup(A) + \Pup(B)$.
    \item \label{item:p4} If $A \cap B = \emptyset$, then $\Plow(A \cup B) \geq \Plow(A) + \Plow(B)$.
\end{enumerate}
However, these properties are not sufficient for coherence.

In the precise case, a probability $p$ and an expectation functional $\mathbb{E}_p[\cdot]$ are in a one-to-one correspondence. We remark that in the imprecise case, the generalized expectation functional, called \textit{coherent lower (upper) prevision}, turns out to be more general,\footnote{Slightly more general even are \textit{sets of desirable gambles} \citep{augustin2014introduction}, which can handle conditioning on events of lower probability zero.} and is therefore the focus of study in the literature (\eg \citep{introtoiplowerprev}): a coherent lower prevision is in a one-to-one correspondence to a closed, convex set of probabilities. For a coherent lower probability, the relation is not one-to-one \citep[Section 2.7.3]{walley1991statistical}. Therefore we will mostly work directly with sets of probabilities. 
\begin{definition}
\label{def:coherentlowerprob}
    We say that functional $\Rlow : \linfty \to \mR$ is a \textit{coherent lower prevision} if
\[
\Rlow(X) = \inf\{\mathbb{E}_{\mu}[X]: \mu \in \mM\}, \quad X \in \linfty,
\]
for some set $\mM$ of probability measures on $\Omega$. Similarly, $\Rup(X)=-\Rlow(-X)$ is the conjugate coherent upper prevision. 
\end{definition} 
Due to the conjugacy, it suffices to consider either $\Rlow$ or $\Rup$ (respectively, $\Plow$ or $\Pup$). The set $\mM$ is unique up to closed convex hull. 
We may sometimes write $\Rlow_{\mM}$ and $\Rlow_{\mM}$ to make the dependence on the set of measures explicit.
We remark that a coherent upper prevision $\Rup$ is equivalently a \textit{coherent risk measure} in the sense of \citet{artzner1999coherent}.

\section{Non-Stationary, Locally Precise Data Models}
\label{sec:nonstationarymodels}
Since we work with a finite possibility set $\Omega$, we can identify each probability measure $p$ on $\Omega$ with a point in the simplex $\Delta^k \coloneqq \{(p_1,p_2,..,p_k) : \sum_{j=1}^k p_j = 1\}$; it induces a linear expectation functional $\mathbb{E}_p[\cdot]$ in the familiar way. We abuse notation and write $p(A)$ when $p \in \Delta^k$, $A \subseteq \Omega$. Due to the finite dimensionality, we assume the topology induced by the Euclidean metric on $\Delta^k$ without loss of generality (see \citep[Appendix B.1]{frohlich2024strictly} for details).
Fix a set $\emptyset \neq \mM \subseteq \Delta^k$.
    For any specified sequence of probability measures $m^\infty=(m_1,m_2,..) \in \mM^\infty$, $m_i \in \mM$, we can define a generator $G$ so that
\[
\lambda\left\{z \in \mathcal{Z} : (\omega_1,\omega_2,..,\omega_n) \subset G(z) \right\} \coloneqq \prod_{i=1}^n  m_i(\{\omega_i\}), \quad \forall n \in \mathbb{N}. 
\]
for some $W_i : [0,1) \rightarrow \Omega$, $i=1..\infty$ and $G : [0,1) \to \Omega^\infty$ with $G(z) \coloneqq (W_1(z),W_2(z),..)$. 
In this model, we keep probabilistic independence but allow the probability to vary ``over time'' (when conceiving of the index $i$ as time), under the constraint that the probabilities belong to $\mathcal{M}$.

It is clear that such a model can exhibit local regularity when $\mM$ is not a singleton. 
More intricate is the question of the aggregate (ir)regularity of such a model. We illustrate this first with a simple example of two coins, taken from \citet{cozmanconvex}.
As a prerequisite, we define the relative frequencies of a data sequence $\omega^\infty=(\omega_1,\omega_2,..)$ as:
\[
r_n^{\omega^\infty} \coloneqq \left(\frac{1}{n} \sum_{i=1}^n \chi_{\{\omega^1\}}(\omega_i), \frac{1}{n} \sum_{i=1}^n \chi_{\{\omega^2\}}(\omega_i), .., \frac{1}{n} \sum_{i=1}^n \chi_{\{\omega^k\}}(\omega_i)\right)^\intercal \in \Delta^k, \; n \in \mathbb{N},
\]
where we stack relative frequencies of each elementary event $\omega^j$, $j=1..k$, in a vector; similarly for the finite case. We write $r_n^{\omega^\infty}(A)$ for $\mathbb{E}_{r_n^{\omega^\infty}}[\chi_A]$, $A \subseteq \Omega$, \ie the relative frequency of the event $A$.

The next examples show that both aggregate regularity and aggregate irregularity can result from such a model \citep{cozmanconvex}.
\begin{example}\normalfont
    Consider two coins, associated with probabilities $1/3$ and $2/3$, respectively, for the event ``heads''. Construct a sequence of probabilities by choosing them in an alternating fashion. Under this independent product, almost surely, the limiting relative frequency of the event ``heads'' is $1/2$. That is, we obtain a stable aggregate, even though there is a hidden heterogeneity, which would be revealed by the selection rule $S(i)=1$ if $i$ is odd, $0$ otherwise and similarly for ``even''.
\end{example}

\begin{example}\label{ex:weirdcoin} \normalfont 
Using the coins as probabilities, construct the sequence as follows. For the $i$-th flip, if the second most significant bit of the binary expansion of $i$ is $0$, take the coin with probability $1/3$ for ``heads''. Otherwise, take the coin with probability $2/3$ for ``heads''. As a result, almost surely there exists no limiting relative frequency for the event ``heads''.
\end{example}

When the limit of relative frequencies does not exist, the natural generalization to consider is a set of cluster points (see \citep[Section 2.3]{frohlich2024strictly} for an argument). Recall that $c \in \mR^k$ is a cluster point of a sequence $a^\infty=(a_1,a_2,..) \in (\mR^k)^\infty$ if
\[
\forall \varepsilon>0 : \forall n_0 \in \mathbb{N}: \exists n \geq n_0: a_n \in B_\varepsilon(c), \quad B_\varepsilon(c) \coloneqq \{x \in \mR^k : d(x,c)<\varepsilon\},
\]
where $d(\cdot,\cdot)$ is the Euclidean metric. In Example~\ref{ex:weirdcoin}, the set of cluster points of relative frequencies for the event ``heads'' is the interval $[4/9,1/2]$. When a topological space is sequentially compact, which $\Delta^k$ with the Euclidean topology is, then a point is a cluster point of a sequence if and only if it is the limit of a convergent subsequence.

Indeed, the focus on relative frequencies of elementary events is justified due to the following.
\begin{proposition}[{\citep[Theorem 1]{ivanenko2017expected},\citep[Proposition 2.1]{frohlich2024strictly}}]
\label{prop:cpsofrelativefrequenciesandgambles}
    Let $\omega^\infty=(\omega_1,\omega_2,..)$ a data sequence and $\ell \in \linfty$ any gamble. Then it holds:
    \[
    \CP\left(n \mapsto \frac{1}{n} \sum_{i=1}^n \ell(\omega_i)\right) = \left\{\mathbb{E}_p[\ell] : p \in \CP(n \mapsto r_n^{\omega^\infty}) \right\}
    \]
\end{proposition}
Here, the operator $\CP(\cdot)$ gives the cluster points of a sequence which is $\mR$-valued (for cluster points of gamble averages) or $\mR^k$-valued (for cluster points of relative frequencies) with respect to the appropriate Euclidean metric in both cases. If $\Omega$ is infinite, a similar statement would hold true when using the weak* topology on the dual space of $\linfty$ but then sequential compactness does not hold (see \citep[Appendix D4]{walley1991statistical}).
In words, knowing the cluster points of relative frequencies suffices to obtain the cluster points of gamble averages for any gamble $\ell \in \linfty$. When a decision maker is interested in long-run loss but a limit does not exist, considering the cluster points provides a natural generalization that may be used for aggregate decision-making.

Our main theorem characterizes exactly which aggregate (ir)regularities can result by choosing from a set of measures to form an infinite independent product.

\begin{theorem}
\label{theorem:maintheoremffivanenko}
    Let $\mathcal{Z}=[0,1)$ and $\emptyset \neq \mathcal{M} \subseteq \Delta^k$. Then there exists a sequence of probability measures $m_1,m_2,.. \in \mM$ with independent product
    \begin{equation}
        \label{eq:independentproduct}
        \lambda\left\{z \in \mathcal{Z} : (\omega_1,\omega_2,..,\omega_n) \subset G(z)\}\right) \coloneqq \prod_{i=1}^n m_i(\{\omega_i\}), \quad \forall n \in \mathbb{N},
    \end{equation}
 so that it holds
    \begin{equation}
    \label{eq:divevent}
    \lambda\left\{z \in \mathcal{Z}:  \CP\left(n \mapsto r_n^{G(z)}\right) =\CP\left(n \mapsto \frac{1}{n} \sumin m_i \right) = \mathcal{N}\right\} = 1
    \end{equation}
if and only if $\emptyset \neq \mathcal{N} \subseteq \cobar(\mathcal{M})$ is a closed connected subset of the closed convex hull of $\mM$. 

We can in addition demand that any $m \in \mM$ is a cluster point of the sequence of measures $m_i$ (with respect to the Euclidean metric).
\end{theorem}
For the proof see Section~\ref{sec:proofofmaintheorem}. The additional demand that any $m \in \mM$ is a cluster point of the $m_i$ sequence only strengthens the statement, meaning intuitively even if any $m \in \mM$ approximately appears infinitely often in the $m_i$ sequence, the set of cluster points of relative frequencies can still be made to coincide with the potentially smaller $\mathcal{N}$.

\begin{definition}
    For brevity, we define $\CM$ as the set of closed connected subsets of the closed convex hull of some $\emptyset \neq \mM \subseteq \Delta^k$. Hence $\mathcal{N} \in \CM$ in the above statement.
\end{definition}
Theorem~\ref{theorem:maintheoremffivanenko} can be viewed as constructing a data model with a typicality notion, which in this case is indeed ``probabilistic'' due to the Lebesgue measure $\lambda$. We define a data model from the converse direction.
\begin{definition}
\label{def:preciseindependentdatamodel}
    Let $\mathcal{Z}=[0,1)$ and $m^\infty=(m_1,m_2,..) \in (\Delta^k)^\infty$. Set 
    \[
    \mM\coloneqq \{m_1,m_2,..\}\;\; \text{ and } \;\; \mathcal{N} \coloneqq \CP\left(n \mapsto \frac{1}{n} \sumin m_i\right) \in \CM.
    \]
    Consider the resulting independent product as in Theorem~\ref{theorem:maintheoremffivanenko}, specified by $G(z)$. 
    We call it the non-stationary, locally precise (\NSLP) data model with $(\mathcal{M},\mathcal{N})$-imprecision and probabilities $m^\infty$. 
    The Lebesgue measure $\lambda$ acts as a typicality notion; formally, $T_a(Z) \coloneqq \lambda(\mathcal{Z} \setminus Z)$, $Z \in \mathcal{B}([0,1))$, is an absolute typicality measure.
\end{definition}
Here, $\mM$ reflects local regularity, whereas $\mathcal{N}$ embodies aggregate irregularity. 
For full rigour we should write ``potentially non-stationary''; note that the label ``non-stationary'' is misleading when $m^\infty$ is an infinite \iid product; of course here we are hardly interested in the \iid models as extreme cases.
Note that we do not make further assumptions on $m^\infty$ here. For instance, a probability measure $m$ need not be a cluster point of $m^\infty$, implying that it cannot be ``revealed'' asymptotically by a selection rule.
As an independent product of precise probability measures, but with potential non-stationarity, we believe this data model has many applications in machine learning and statistics (see Section~\ref{sec:applications}).
This view is only a mild departure from the familiar \iid model, since at each ``timestep'' $i \in \mathbb{N}$, data is still modelled as drawn from a precise probability, but the probability varies over ``time''. Hence imprecision cannot  be considered a local characteristic of the data model. In other words, the imprecision manifests itself through non-stationarity in the forms of aggregate (ir)regularity and local (ir)regularity.

Under this data model, parameterized by a sequence of measures, diverging relative frequencies with cluster points $\mathcal{N}$ are typical. Observe that to obtain such typical divergence of relative frequencies, the sequence $m^\infty$ cannot itself converge, since then $\mathcal{N}$ must be a singleton that coincides with this limit: in other words, perfect local irregularity (as in the \iid model) cannot yield aggregate irregularity.

As a trivial corollary of Theorem~\ref{theorem:maintheoremffivanenko}, we can characterize the limiting behaviour of subsequences extracted by selection rules. Denote by $\omega_S^\infty$ the subsequence of $\omega^\infty$ extracted by some selection rule $S$ (keeping outcomes at those indices for which $S(i)=1$).
\begin{corollary}
\label{cor:cpsofselectionrules}
    Let $m^\infty=(m_1,m_2,..) \in (\Delta^k)^\infty$ as in Definition~\ref{def:preciseindependentdatamodel} and the corresponding non-stationary, locally precise data model, and $S: \mathbb{N} \to \{0,1\}$ a selection rule.
    Then there exists a unique $\mathcal{N} \in \CM$ so that
      \[
    \lambda\left\{z \in \mathcal{Z}: \omega^\infty \coloneqq G(z) : \CP(n \mapsto r_{n}^{\omega_S^\infty}) = \mathcal{N}\right\} = 1.
    \]
\end{corollary}
The implication: even if the manifestation we are interested in is local (ir)regularity, the result on aggregate (ir)regularity is useful since it tells us what behaviour we can expect along selection rules: asymptotically, relative frequencies along a selection rule never end up outside the closed convex hull of $\mM$. Intuitively, any hidden heterogeneity is confined to lie therein. Since a selection rule induces an aggregate in its own right, a subsequence, it can itself exhibit aggregate (ir)regularity --- indeed it is even possible to have overall aggregate regularity, but that a selection rule yields aggregate irregularity.\footnote{Consider the infinite binary sequence $\omega^\infty=(0,1,0,1,..)$, with aggregate regularity. Acknowledging the existence of binary sequences with diverging relative frequencies implies that we can find a selection rule that yields aggregate irregularity.} Furthermore, if some $m \in \mM$ is a cluster point of the sequence of measures, there exists a selection rule which ``reveals'' it.

\begin{proposition}
\label{prop:uncoverhiddenheterogeneity}
    Let $m^\infty=(m_1,m_2,..) \in (\Delta^k)^\infty$ so that $m \in \CP(i \mapsto m_i)$. Then there exists a selection rule $S: \mathbb{N} \to \{0,1\}$ so that
\[
 \lambda\left\{z \in \mathcal{Z}: \omega^\infty \coloneqq G(z) : \CP(n \mapsto r_{n}^{\omega_S^\infty}) = m \right\} = 1.
\]
\end{proposition}
\begin{proof}
    Take any sequence $\varepsilon^\infty=(\varepsilon_1,..) \downarrow 0$. Since $m \in \CP(i \mapsto m_i)$, we can, to any $\varepsilon_i$ find some $i^*$ so that $d(m_{i^*},m)<\varepsilon_i$ and that the sequence of the $i^*$ is strictly increasing. Choose these $i^*$ to construct a selection rule, \ie $S(j)=1$ if $j=i^*$ for some $\varepsilon_i$, $0$ otherwise. That is, this selection rule selects a subsequence of measures which converges to $m$, but then also the sequence of averages of those measures converges to $m$. Applying Theorem~\ref{theorem:maintheoremffivanenko} on this subsequence then gives the result.
\end{proof}
Observe that in general, if $p \in \mathcal{N}$ is \textit{not} in $\CP(i \mapsto m_i)$ but in $\CP(n \mapsto \frac{1}{n} \sumin m_i)$, then clearly there need not exist a selection rule that reveals it as local irregularity in the above sense.

As a special case of Theorem~\ref{theorem:maintheoremffivanenko}, we obtain Theorem~3.3 in \citep{frohlich2024strictly}.
\begin{corollary}
    Take $\mM=\{e_1,..,e_k\}$ where $e_i=(0,..,1,..0)$ with the $1$ at the $i$-th position, that is, the vertices of the simplex. Then Theorem~\ref{theorem:maintheoremffivanenko} asserts that if and only if $\mathcal{N} \subseteq \Delta^k$ is closed and connected, there exists a sequence $\omega^\infty$ so that $\CP\left(n \mapsto r_n^{\omega^\infty}\right) = \mathcal{N}$.
\end{corollary}
In contrast to us, \citet{frohlich2024strictly} have worked in a \textit{strictly frequentist setting}, meaning that they have taken the sequence $\omega^\infty$ as the primitive entity, and defined a coherent lower (upper) prevision based on the set of cluster points of $r_n^{\omega^\infty}$. In this way, they follow up on the work of \cite{ivanenkobook}. In their approach, data models do not appear. Put differently, process randomness and typicality have no place in their framework, only outcome randomness.

\subsection{A Comparison to the Fierens-Fine Model}
\label{sec:fierensfinecomparison}
It is also pertinent to compare the data model of Definition~\ref{def:preciseindependentdatamodel} to the works of \citet{fierens2001towards,fierens2003towards,fierens2003towardsphdthesis,rego2005estimation,fierens2009extension,fierens2009frequentist}. In the early paper \citep{fierens2001towards} the authors focus on the manifestation of aggregate (ir)regularity, as they study ``persistent oscillations of relative frequencies'', following up on a previous line of work which we discuss in Section~\ref{sec:comparisontofine}. However the paper \citep{fierens2003towards} marks a conceptual turning point:\footnote{\citet{fierens2003towards} write ``We subsequently judged this
approach to be inadequate, in part after considering the performance of martingale
betting systems on such time series as advocated by the then newly-published [book by]
Shafer and Vovk [..].''; they refer here to the game-theoretic account of probability put forward by \citet{shafervovk2001probabilityfinance}.}
from there on the work focuses on the manifestation local (ir)regularity. Why? Our understanding of their motivation for disregarding aggregate (ir)regularity is that they wanted ``to ensure that there was no further exploitable structure in the time evolution'' \citep{fierens2003towards}. Yet as \citet[p.\@\xspace 43]{fierens2003towardsphdthesis}  points out in the language of game-theoretic probability, (in our reading) a skeptic who is allowed to take up betting commitments from a decision maker whose estimation method is focused on aggregate (ir)regularity can exploit local (ir)regularity to get infinitely rich. This can be illustrated with Example~\ref{ex:weirdcoin}: here, the underlying probabilities are the set $\{1/3,2/3\}$, and the cluster points of relative frequencies are the interval $[4/9,1/2]$. If a skeptic knows the pattern for which coin is chosen based on $i \in \mathbb{N}$, and the decision maker offers fair bets based on the interval $[4/9,1/2]$,\footnote{Such a bet is ``fair'' when it consists in playing a \textit{marginally desirable gamble} \citep{walley1991statistical}. For example, the gamble $X-\Pup(X)$ and $\Plow(X)-X$ are marginally desirable for a decision maker with coherent upper prevision $\Pup$. In the precise case, this is equivalent to $\mathbb{E}_p[X-\mathbb{E}_p[X]]=0$, \ie $\mathbb{E}_p$ determines betting rates.}, skeptic can get infinitely rich when gambling infinitely often since $[4/9,1/2] \subsetneq [1/3,2/3]$. We think this is hardly surprising, and serves to highlight that whether a focus on aggregate (ir)regularity or local (ir)regularity is warranted depends on the decision-making context, see Section~\ref{sec:applications}.
In \citep{fierens2003towards} and subsequent papers,
the authors then consider a data model similar to the \NSLP~model, and compute relative frequencies along selection rules. Their model is more expressive, however, since they allow the choice of probabilities $p_i$ to depend on previously observed outcomes. In this paper, for simplicity we focus on the case where the sequence of probabilities is fixed a priori. We believe that Theorem~\ref{prop:mainproptheoretical} can be extended to their setting, but leave the details to the interested reader.\footnote{To this end, consider the following variant of the strong law of large numbers given by Loève (1978, as cited in \citep{papamarcou1991unstable}): For an increasing sequence of $\sigma$-algebras $\mathcal{B}_i$ and a uniformly bounded and adapted sequence of measurable functions $Y_i$, it holds that $\limn \frac{1}{n} \sumin \left(Y_{i+1} -\mathbb{E}[Y_{i+1}|\mathcal{B}_i]\right)=0$ \textit{almost surely}.}
We also note that the earlier paper by \citet{cozmanconvex} represents a ``bridge'' between the two manifestations, since they consider both subsequences of the relative frequency sequence as well as subsequences of the data sequence.

Starting with the non-stationary, locally precise data model, more exotic data models can be obtained, where the imprecision is already inherent to each ``timestep'' itself, meaning it is a local feature. As we will see, imprecision can therefore manifest even under \textit{stationarity}. Investigating such models is also of interest since it explicates the relation of Theorem~\ref{theorem:maintheoremffivanenko} and the present data model to the works of \citet{walley1982towards} and subsequent papers \citep{kumar1985stationary,grize1987continuous,fine1988lower,papamarcou1991stationarity,papamarcou1991unstable,sadrolhefazi1994finite}. One of our goals is to provide an accessible introduction to this line of work, therefore we now consider such more general models.

\section{Stationary, Locally Imprecise Data Models}
\label{sec:stationarymodels}
Previously we have fixed a specific sequence of probability measures from $\mM$. 
Like \citet[Section 5]{kumar1985stationary}, we now construct data models which subsume whole sets of sequences of measures. The resulting models are highly related to those of \citet{walley1982towards}. 
Let $\emptyset \neq \mM \subseteq \Delta^k$ be a set of probability measures. We write $\Plow_{\mM}(A) \coloneqq \inf\{m(A) : m \in \mM\}$, $A \subseteq \Omega$, and recall that we denote by $\mathcal{M}^\infty$ the set of sequences of probability measures from $\mM$.

Given a subset $\emptyset \neq \tilde{\mM}^\infty \subseteq \mM^\infty$ we construct a data model based on independent products. Define the infinite product space of $\Omega$ and the associated product $\sigma$-algebra as:
\[
\Omega^\infty \coloneqq \bigtimes_{i=0}^\infty \Omega, \quad \bigtimes_{i=0}^\infty  2^\Omega \coloneqq \sigma\left(\bigcup_{i=1}^\infty \sigma(\pi_i)\right), 
\]
where the projections are $\pi_i : \Omega^\infty \to \Omega$, $\pi_i(\omega^\infty=(\omega_1,..)) \coloneqq \omega_i$.
This product $\sigma$-algebra coincides with the Borel $\sigma$-algebra $\infprodsigma$ when assuming the product topology induced by the discrete topology on $\Omega$. For terse notation, we therefore use $\infprodsigma$.
A set $A^\infty \in \infprodsigma$ is a set of infinite sequences; but note that $\infprodsigma \subsetneq 2^{\Omega^\infty}$. Fix any sequence of probability measures $m^\infty = (m_1,m_2,..) \in \tilde{\mM}^\infty$. 
Ionescu-Tulcea's theorem \citep[Section 14.3]{klenke2013probability} guarantees the existence of a unique measure $\operatorname{IT}(m^\infty)$ on $(\Omega^\infty,  \infprodsigma)$ which satisfies
 \begin{equation}
      \label{eq:rectangles}
\operatorname{IT}(m^\infty)\left(A_1 \times .. \times A_n \times \bigtimes_{i=n+1}^\infty \Omega \right) \coloneqq \prod_{i=1}^n m_i(A_i), \quad \forall n \in \mathbb{N}_0.
 \end{equation} 
In particular, this implies
 \[
\operatorname{IT}(m^\infty)\left(\{\omega^\infty = (\omega_1, \omega_2, .., \omega_n) \}\right) \coloneqq \prod_{i=1}^n m_i(\{\omega_i\}), \quad \forall n \in \mathbb{N}_0.
\] 
For each such sequence $m^\infty \in \tilde{\mM}^\infty$, we obtain a measure $\operatorname{IT}(m^\infty)$ on $(\Omega^\infty, \infprodsigma)$ from Ionescu-Tulcea's theorem. Since these measures are all defined on the same probability space, we can define a coherent upper probability from them via an envelope representation.

Depending on how we constrain the set $\tilde{\mM}^\infty$, we can obtain different data models. For instance, based on a set of probability measures $\mM$, we could consider \textit{all} sequences of measures from $\mM$, yielding the following.
 \begin{definition}
   Let $\emptyset \neq \mM \subseteq \Delta^k$. 
   We define the stationary, locally imprecise (\SLI) data model with $\mM$-imprecision via a coherent lower probability:
   \[
\Plow_{\mM}^\infty(A^\infty) \coloneqq \inf\left\{\operatorname{IT}(m^\infty)(A^\infty) : m^\infty \in  \mM^\infty\right\}, \quad \forall A^\infty \in \infprodsigma,
\]
with conjugate coherent upper probability $\Pup_{\mM}^\infty(A^\infty) = 1-\Plow_{\mM}^\infty((A^\infty)^c)$.
\end{definition}
This model is essentially (see the discussion below) equivalent to the \iid model of \citet{walley1982towards}.
We have here expressed a notion of typicality on $\infprodsigma$ through a coherent lower probability: intuitively, $\Plow_{\mM}^\infty(A^\infty)=1$ means that $A^\infty$ is a typical set, since this implies that all measures in the envelope assign measure $1$ to the set. We can formalize this in the framework of \citet{galvan2006bohmian}.
\begin{proposition}
\label{prop:galvantypicality}
Let $\Pup$ be a coherent upper probability, defined on some possibility set $\Lambda$ and set system $\mathcal{A}$, which is closed under complementation and $\emptyset \in \mathcal{A}$. Then $d(A,B)\coloneqq \Pup(A \triangle B)$ satisfies the axioms of a typicality distance as in \citep{galvan2006bohmian}. Thus $(\Omega,\mathcal{A},d)$ forms a typicality space \citep{galvan2006bohmian}.The induced absolute typicality operator is $T_a(A) = 1-\Plow(A)$, where $\Plow$ is conjugate to $\Pup$.
\end{proposition}
\begin{corollary}
    If $d(A^\infty,B^\infty) \coloneqq \Pup_{\mM}^\infty(A^\infty \triangle B^\infty)$, then $(\Omega^\infty,\infprodsigma,d)$ forms a typicality space.
\end{corollary}
For the definitions of the involved concepts and the proof, see Appendix~\ref{app:galvantypicalityproof}
Thus, $T_a(A^\infty)=0 \Leftrightarrow \Plow_{\mM}^\infty(A^\infty)=1$ means that $A^\infty$ is a typical set.
To better fit the notion of an abstract random seed, we could have also defined it on $[0,1)$ instead.\footnote{One possibility would be embedding $\Omega$ into $\mR$, so that the $m_i$ would share the same finite support, and then exploit that $[0,1)$ and $\mR^\infty$ are Borel-isomorph.}

In the above model, we have allowed arbitrary sequences of probability measures from $\mM$. 
Another possibility is to constrain the set by a specific set of cluster points. Consider the subset 
\[\mM_{\shortrightarrow \mathcal{N}}^\infty \coloneqq \left\{m^\infty = (m_1,m_2,..) \in \mM^\infty : \CP\left(n \mapsto \frac{1}{n} \sum_{i=1}^n m_i \right) = \mathcal{N}\right\},
\] for some $\mathcal{N} \in \CM$, \ie a closed connected subset of the closed convex hull of $\mM$. 
In this set, all sequences of probability measures from $\mM$ are collected which yield a specific set $\mathcal{N}$ as cluster points. This leads to the following more constrained model.

 \begin{definition}
 \label{def:datamodelwithmnimprecision}
  Let $\emptyset \neq \mM \subseteq \Delta^k$ and $\mathcal{N} \in \CM$. 
  We define the stationary, locally imprecise data model with $(\mM,\mathcal{N})$-imprecision through a notion of typicality expressed by a coherent lower probability:
   \[
\Plow_{\MtoNset}^\infty(A^\infty) \coloneqq \inf\left\{\operatorname{IT}(m^\infty)(A^\infty) : m^\infty \in  \mM_{\shortrightarrow \mathcal{N}}^\infty\right\}, \quad \forall A^\infty \in \infprodsigma.
\]
\end{definition}

Clearly, cluster points of relative frequencies will typically coincide with $\mathcal{N}$ in this data model.
\begin{proposition}
\label{prop:fullsupport}
    Consider the event $\DtoN \coloneqq \{\omega^\infty: \CP(n \mapsto r_n^{\omega^\infty}) = \mathcal{N}\} \in \infprodsigma$. Then $\Plow_{\MtoNset}^\infty(\DtoN) = \Pup_{\MtoNset}^\infty(\DtoN)= 1$, that is, $\DtoN$ is a typical set.
\end{proposition}
\begin{proof} First observe that indeed $\DtoN \in \infprodsigma$. Next, by construction, if $m^\infty \in \MtoNset$, then $\operatorname{IT}(m^\infty)(\DtoN)=1$. But then taking the infimum gives $\Plow_{\MtoNset}^\infty(\DtoN)=1$.
\end{proof}

\subsection[A Comparison to Walley \& Fine's (1982) and Fine et al.'s line of work]{A Comparison to \citeauthor{walley1982towards}'s \citeyearpar{walley1982towards} and Fine et al.'s line of work}
\label{sec:comparisontofine}
It is instructive to compare our models to the seminal work of \citet{walley1982towards} and the subsequent line of work by Fine and collaborators (starting with \citep{kumar1985stationary}). We begin with a comparison to \cite{walley1982towards}, which yields important insight into the behaviour of these models. We consider their imprecise \iid model, which is also what they focus on. 
In contrast to our coherent lower probability, theirs is defined on the set system \[
\mathcal{A}^\infty \coloneqq \bigcup_{n \in \mathbb{N}} \mathcal{A}^n, \quad \text{ where } \mathcal{A}^n \coloneqq \sigma\left(\bigcup_{i=1}^n \sigma(\pi_i)\right).
\]
Recall that the projections are $\pi_i : \Omega^\infty \to \Omega$, $\pi_i(\omega^\infty=(\omega_1,..)) \coloneqq \omega_i$. 
Note that $\mathcal{A}^\infty$ is an \textit{algebra} (also called \textit{field}), but not a $\sigma$-algebra. 
Since $\sigma(\mathcal{A}^\infty)=\infprodsigma$, their set system is a strict subset of ours.
Our model partially agrees with the \iid model of \citet{walley1982towards} in the following sense. We denote by $\PlowWF : \mathcal{A}^\infty \to [0,1]$ their \iid model, which is also parameterized by a set of measures $\mM$. Intuitively, this model corresponds to $\Plow_{\mM}^\infty$, since it is built from independent lower probabilities, each of which is the envelope of $\mM$; for the detailed construction, we refer the reader to \citet{walley1982towards}. However, for events which depend only on a finite number of outcomes, the model agrees with both $\Plow_{\MtoNset}^\infty$ and $\Plow_{\mM}^\infty$

\begin{proposition}
\label{prop:modelonrectangles}
    Let $A^\infty \coloneqq A_1 \times .. \times A_n \times \bigtimes_{i=n+1}^\infty \Omega$ be a rectangle in $\infprodsigma$, $A_i \in 2^\Omega$. It holds that
    \[
    \Plow_{\MtoNset}^\infty\left(A^\infty \right) = \Plow^\infty_{\mM}(A^\infty) = \PlowWF(A^\infty) = \inf\left\{\prod_{i=1}^n m_i(A_i) : m_i \in \mathcal{M}\right\} , \quad \forall n \in \mathbb{N}_0.
    \]
\end{proposition}
\begin{proof}
   The part of the statement regarding $\PlowWF(A^\infty)$ follows from its definition, see \citet{walley1982towards}. The statement regarding $\Plow_{\mM}^\infty$ is clear from its definition. For $\Plow_{\MtoNset}^\infty$, we know that
    \[
    \Plow_{\MtoNset}^\infty\left(A_1 \times .. \times A_n \times \bigtimes_{i=n+1}^\infty \Omega\right) = \inf\left\{\prod_{i=1}^n m_i(A_i) : m^\infty=(m_1,m_2,..) \in \MtoNset\right\}, \quad \forall n \in \mathbb{N}_0.
    \]
    But since the cluster points of $m^\infty$ remain unchanged when prepending an arbitrary prefix to the sequence of measures, the statement follows obviously: any chosen prefix $(m_1,..,m_n) \in \mM^n$ can be extended to lie in $\MtoNset$.
\end{proof}
The algebra that is generated by sets of the above form, the rectangles, is $\mathcal{A}^\infty$, which implies the following.
 \begin{corollary}
     Let $A^\infty \in \mathcal{A}^\infty \subsetneq \infprodsigma$. Then 
     \[
     \Plow_{\MtoNset}^\infty(A^\infty) = \Plow^\infty_{\mM}(A^\infty) = \Plow^\infty_{\mM,\operatorname{WF}}(A^\infty).
     \]
 \end{corollary}
 \begin{proof}
     Rectangles as above form a $\pi$-system, \ie a set system closed under finite intersections, and hence the measures in the envelope, which agree pairwise on this $\pi$-system, extend uniquely to algebra generated by it (according to the $\pi$-$\lambda$-theorem), which is $\mathcal{A}^\infty$. By pairwise agreement, we mean that for any measure in the envelope of $\Plow_{\mM}$ ($\Plow_{\MtoNset}$, respectively), we can find a measure in the envelope of $\PlowWF$ so that these agree on the $\pi$-system, and conversely.
 \end{proof}
In contrast, for events $A^\infty \in \infprodsigma \setminus \mathcal{A}^\infty$, the assessments of the models can be starkly different. 
Define the tail $\sigma$-algebra induced by the projections as
\[
\mathcal{T} \coloneqq \bigcap_{n \in \mathbb{N}} \sigma\left(\pi_{n+1}, \pi_{n+2}, ..\right), \quad \mathcal{T} \subsetneq \infprodsigma.
\]
Intuitively, $\mathcal{T}$ is the sub-$\sigma$-algebra of $\infprodsigma$ which does not depend on the outcomes of finitely many trials. It is clear that $A^\infty \in \mathcal{T} \Longrightarrow A^\infty \notin \mathcal{A}^\infty$; however note that $\infprodsigma \setminus \mathcal{T} \neq \mathcal{A}^\infty$, implying that $\mathcal{T} \subsetneq \infprodsigma \setminus \mathcal{A}^\infty$. For illustration, consider the event $\{\omega^\infty: \omega^\infty = \omega^\infty_0\}$, asserting that the sequence is exactly a specified sequence $\omega^\infty_0$. This event is neither in $\mathcal{T}$ nor in $\mathcal{A}^\infty$.

In a sense, $\Plow_{\mM}^\infty$ can be viewed as an extension of the Walley-Fine \iid model to all of $\infprodsigma$, so that the event of divergence is a measurable tail event: $\DtoN \in \mathcal{T} \subsetneq \infprodsigma$.
\begin{proposition} Let $|\mM|>1$. It holds that:
   \begin{enumerate}[label=(\theproposition.\arabic*),leftmargin=1.4cm]
  \item \label{item:fullsupport} $\Plow_{\MtoNset}^\infty(\DtoN)=\Pup_{\MtoNset}^\infty(\DtoN)=1$.
  \item \label{item:weaksupport} $\Plow_{\mM}^\infty(\DtoN)=0$ but $\Pup_{\mM}^\infty(\DtoN)=1$.
   \end{enumerate}
\end{proposition}
\begin{proof}
Property~\ref{item:fullsupport} is just a restatement of Proposition~\ref{prop:fullsupport}.
To see that $\Plow_{\mM}^\infty(\DtoN)=0$, assume first that $\mathcal{N} \subsetneq \cobar(\mM)$. Then consider the sequence $m^\infty=(m_0,m_0,..)$ for some $m_0 \in \mM \cap \mathcal{N}^c$. Then $IT(m^\infty)(\DtoN)=0$. On the other hand, if $\mathcal{N} = \cobar(\mM)$, then take any $m_0 \in \mM$, yielding $IT(m^\infty)(\DtoN)=0$.
Finally, $\Pup_{\mM}^\infty(\DtoN)=1$ follows from Theorem~\ref{theorem:maintheoremffivanenko}. 
\end{proof}
We emphasize again that throughout this section we assume that $\mathcal{N}$ is any set that satisfies the appropriate relation to $\mM$, \ie $\mathcal{N} \in \CM$.
Intuitively, the Walley-Fine \iid model is completely silent about the set of cluster points, about convergence or divergence. 
Note that in general $[\Plow_{\MtoNset}(A^\infty), \Pup_{\MtoNset}(A^\infty)] \subseteq [\Plow_{\mM}(A^\infty), \Pup_{\mM}(A^\infty)]$ since $\MtoNset \subsetneq \mM^\infty$ for any $\mathcal{N} \in \CM$.

\begin{example}
    Take $\mM=\{e_1,..,e_k\}$ where $e_i=(0,..,1,..0)$ with the $1$ at the $i$-th position. Then a suitable $\mathcal{N} \in \CM$ is any closed connected subset of $\Delta^k$. 
    Each infinite sequence for which cluster points of relative frequencies coincide with $\mathcal{N}$ then corresponds to an infinite product measures in $\MtoNset$. Hence $\DtoN$ is a typical set. Similarly, $\Plow_{\mM}$ is a fully vacuous model which is silent on all ``nontrivial'' events, since it is the envelope of all $\Omega^\infty$, as represented by Dirac-measures.
\end{example}

While $\Plow_{\mM}^\infty$ is a reasonable extension of the model to $\infprodsigma$, the original proposal by \citet{walley1982towards} made the architectural choice to work on $\mathcal{A}^\infty$ as opposed to $\infprodsigma$. One reason is that the scope of \citet{walley1982towards} also includes non-\iid models in principle (dependent coherent lower probabilities), which are built up from a consistent sequence of coherent lower probabilities on $\mathcal{A}^1,\mathcal{A}^2,..$; but the focus of \citet{walley1982towards} is strongly on the \iid model. The restriction to $\mathcal{A}^\infty$ means that ``strong law''-related statements\footnote{In the terminology of \citep{kumar1985stationary} and meaning roughly \textit{in relation to a generalization of the strong law of large numbers}.} of the above form cannot be made. Instead, \citet{walley1982towards} have studied the event of \textit{apparent convergence} as
\[
C_n(A; k, \varepsilon) \coloneqq \bigcap_{j=k}^n \left\{\omega^\infty:  |r_j^{\omega^\infty}(A) - r_n^{\omega^\infty}(A)| < \varepsilon\right\}.
\]
In words, ``the event that $r_1(A),..,r_n(A)$ apparently converges $(k,\varepsilon)$''. Considering all events simultaneously, they define $C_n(k,\varepsilon)\coloneqq \bigcap_{A \in 2^\Omega} C_n(A;k,\varepsilon)$. Similarly, the event of apparent divergence is $D_n(k,\varepsilon)\coloneqq C_n(k,\varepsilon)^c$. Note that this event is in $\mathcal{A}^\infty$, and thus our model agrees with that of \citet{walley1982towards}; we then state results in terms of $\Plow_{\mM}^\infty$. We emphasize that the nature of the event $D_n(k,\varepsilon) \in \mathcal{A}^n$ is fundamentally different to $\DtoN \in \mathcal{T} \cap (\mathcal{A}^\infty)^c$; intuitively, $D_n(k,\varepsilon)$ is a ``weak law''-related statement. Recall that the weak law of large numbers in the precise \iid case asserts the following:
\[
\forall \varepsilon>0: \limn P^\infty\left(|r_n^{\omega^\infty}(A) - P(A)|>\varepsilon\right)=0, \quad P^\infty=\operatorname{IT}((P,P,..)).
\]

In the imprecise case, the situation is somewhat unsettling. Recall that we write $\Plow_{\mM}(A) \coloneqq \inf\{m(A) : m \in \mM\}$, $A \subseteq \Omega$, for the envelope of $\mM$. Intuitively, this is the marginal of $\Plow_{\mM}^\infty$.
\begin{proposition}[{\citep[Theorem 4.1]{walley1982towards}}] Let $\kfunc: \mathbb{N} \to \mathbb{N}$, $\kfunc(n)\to \infty$ as $n\to \infty$.
    For any $\varepsilon>0$, 
    \begin{enumerate}[label=(\theproposition.\arabic*),leftmargin=1.6cm]
    \item \label{item:lowernosupportfordiv} $    \limn \PlowWF(D_n(\kfunc(n),\varepsilon)) = 0.$
    \item \label{item:uppersupportfordiv} If $|\mM|>1$, $\limsupn \kfunc(n)/n = 0$, and $\exists A \subseteq \Omega$ : $\varepsilon< \Pup_{\mM}(A) - \Plow_{\mM}(A)$, then \[\limn \Pup_{\mM}^\infty(D_n(\kfunc(n),\varepsilon)) = 1.\]
    \end{enumerate}
\end{proposition}
From~\ref{item:lowernosupportfordiv}, it follows that $    \limn \Plow_{\MtoNset}(D_n(\kfunc(n),\varepsilon)) = 0$.
The situation is a curious one: while a ``strong law''-type statement holds, that divergence with a set of cluster points $\mathcal{N}$ is typical under $\Plow_{\MtoNset}^\infty$, a corresponding ``weak law'' fails. To address this, \citet{walley1982towards} coined a substantially weaker notion of typicality, called \textit{asymptotic favorability}, so that $D_n(\kfunc(n),\varepsilon)$ is found to be asymptotically favorable; the statement~\ref{item:uppersupportfordiv} is similar to Theorem 4.1.d in \citet{walley1982towards}, which expresses asymptotic favorability of $D_n(\kfunc(n),\varepsilon)$: a sequence of events $A_n$ is called asymptotically favorable if $\Plow^\infty(A_n^c) / \Plow^\infty(A_n) \to 0$, which implies $\Pup^\infty(A_n) \to 1$. The intuitive reason for why a ``weak law''-type statement must fail in this model is that for any rate of divergence, there exists some other product measure which yields even slower divergence. Hence statements regarding apparent divergence \textit{must} be vacuous (for appropriate $\varepsilon$). In contrast, to anticipate confusion for when we later contextualize our results and compare to existing weak laws for imprecise probabilities in the literature, we give a slighty modified restatement of \citep[Theorem 4.1.a]{walley1982towards}.
\begin{proposition}
\label{prop:wfweaklaw}
    For any $\varepsilon>0$,
    \[
    \limn \Plow_{\mM}^\infty(\{\omega^\infty: \forall A \subseteq \Omega: \Plow_{\mM}(A) - \varepsilon < r_n^{\omega^\infty}(A) < \Pup_{\mM}(A) + \varepsilon \}) = 1.
    \]
\end{proposition}
A similar proposition holds for gamble averages by replacing the lower/upper probability with the lower/upper prevision. Intuitively, this means that while ``weak law'' statements about divergence fail, \ie the model is silent on them, the precise probabilistic weak law\footnote{A version of the weak law for independent but non-\iid measurable functions.} does hold for all of the product measures in the set of measures, and therefore we find that with increasing lower probability, relative frequencies must lie between the lower and upper probability. This should be rather unsurprising. Thus, in this section and in contrast to Section~\ref{sec:llncomparisonsection}, when we speak of ``strong law'' or ``weak law''-type statements, we refer to statements that concern events of divergence or convergence, instead of merely a confinement to the interval determined by the marginal $\Plow_{\mM}$.

The work of \citet{walley1982towards} sparked a long and difficult to access line of work by Fine and collaborators, \citep{kumar1985stationary,grize1987continuous,fine1988lower,papamarcou1991stationarity,papamarcou1991unstable,sadrolhefazi1994finite}, which took the above phenomenon, the lack of typical divergence, as the starting point.  \citet{fine1988lower} provides a relatively accessible overview to the philosophy and approach of this line of work. These works also moved from the setup $(\Omega^\infty,\mathcal{A}^\infty)$ to $(\Omega^\infty,\infprodsigma)$, like ours.
The authors were motivated by the study of physical systems which seemingly display the puzzling combination of stationarity and unstable averages. \citet{kumar1985stationary} define stationarity as follows.

\begin{definition}
     For $k \in \mathbb{N}$, $T^k: \Omega^\infty \to \Omega^\infty$ is the \textit{left-shift-by-$k$} operator defined by \[T^k(\omega^\infty=(\omega_1,\omega_2,..))\coloneqq (\omega_{k+1},\omega_{k+2},..),
     \] which acts on a set $A^\infty \in \infprodsigma$ as $T^k(A^\infty) \coloneqq \{T^k(\omega^\infty) : \omega^\infty \in A^\infty\}$. Conversely, $T^{-k}$ is defined by $T^{-k}(A^\infty) \coloneqq \{\omega^\infty: T^k(\omega^\infty) \in A^\infty\}$, $A^\infty \in \infprodsigma$.
\end{definition}
\begin{definition}
     A set function $\Plow$ on $(\Omega^\infty,\infprodsigma)$ is called  stationary when 
     \[
     \forall A^\infty \in \infprodsigma : \forall k \in \mathbb{N}: \Plow(T^{-k}(A^\infty)) = \Plow(A^\infty).
     \]
\end{definition}
Indeed, the models presented in this section are stationary in this sense (compare Example 5.1 of \citet{kumar1985stationary}).
\begin{proposition}
    $\Plow_{\mM}^\infty$ and $\Plow_{\MtoNset}^\infty$ are stationary.
\end{proposition}
\begin{proof}
We first show stationarity of $\Plow_{\mM}$. 
    Consider an arbitrary $m^\infty=(m_1,m_2..) \in \mM^\infty$. We show that if $\operatorname{IT}(m^\infty)(A^\infty)=a$, there is some other $\tilde{m}^\infty \in \mM^\infty$ so that $\operatorname{IT}(\tilde{m}^\infty)(T^{-k}(A))=a$. From this we then conclude $\Plow_{\mM}(T^{-k}(A)) \leq \Plow_{\mM}(A)$. 
    Take 
    \[
    \tilde{m}^\infty=(\underbrace{m_0,..,m_0}_{k \text{ times}},m_1,m_2..)
    \]
    for arbitrary $m_0 \in \mM$. Consider rectangles of the form 
    \[
    \mathcal{R}_n \coloneqq \left\{ 
    \underbrace{\Omega \times .. \times \Omega}_{k \text{ times}} \times A_1 \times .. \times A_n \times \bigtimes_{k+n+1}^\infty \Omega : A_i \in 2^\Omega \right\}, \quad n \in \mathbb{N}.
    \]
    These rectangles $\{\mathcal{R}_n : n \in \mathbb{N}\}$ generate the $\sigma$-algebra $\mathcal{A}^{k+1:} \coloneqq \sigma\left(\bigcup_{i=k+1}^\infty \sigma(\pi_i)\right)$. It is clear that for these rectangles, $R \mapsto \operatorname{IT}(\tilde{m}^\infty)(R)$ coincides with $R \mapsto m_0(\Omega) \cdot m_0(\Omega) \cdot \operatorname{IT}(T^k(R))$ due to \eqref{eq:rectangles}, and $m_0(\Omega)=1$. Since these rectangles form a $\pi$-system, these measures also coincide on the generated $\sigma$-algebra  $\mathcal{A}^{k+1:}$. Then, noting that $T^{-k}(A^\infty) \in \mathcal{A}^{k+1:}$, we have 
    \[
    \operatorname{IT}(\tilde{m}^\infty)(T^{-k}(A^\infty)) = \operatorname{IT}(m^\infty)(T^k(T^{-k}(A^\infty))),
    \] which since $T^{k}(T^{-k}(A^\infty)=A^\infty$ implies $\operatorname{IT}(\tilde{m}^\infty)(T^{-k}(A^\infty))=\operatorname{IT}(m^\infty)(A^\infty)$; we caution the reader that $T^{-k}(T^k(A^\infty)) \supseteq A^\infty$. 

    For the converse direction we want to show that $\Plow_{\mM}(A^\infty) \leq \Plow_{\mM}(T^{-k}(A))$ for $A^\infty \in \infprodsigma$. For this we exhibit to any $m^\infty=(m_1,m_2,..) \in \mM^\infty$ with $\operatorname{IT}{m^\infty}(T^{-k}(A^\infty))=a$ some $\tilde{m}^\infty \in \mM^\infty$ with $\tilde{m}^\infty(A^\infty)=a$. But here it is clear that $\tilde{m}^\infty \coloneqq (m_{k+1},m_{k+2},..)$ works.

    Finally, for the stationarity of $\Plow_{\MtoNset}$, observe that the logic of the proof works as well with the constraint $m^\infty \in \MtoNset$ (compare the proof of Proposition~\ref{prop:modelonrectangles}).
\end{proof}

We need some preliminaries to state the main result of \cite{kumar1985stationary}. The authors define a lower/upper probability as a conjugate pair $\Plow(A)=1-\Pup(A)$ of set functions, satisfying a list of axioms (implying~\ref{item:p1}-~\ref{item:p4}), which are weaker than the coherence condition of \citet[Section 2.5]{walley1991statistical}. To a lower probability, associate the set
\[
\mM_{\Plow} \coloneqq \{\mu : \mu \in \Delta^k \text{ and } \mu(A) \geq \Plow(A) \;\; \forall A \subseteq \Omega\}.
\]
When this set is nonempty, call $\Plow$ a \textit{dominated} lower probability (equivalent to the \textit{avoiding sure loss} condition in \citet[Section 2.4]{walley1991statistical}); when it is empty, call $\Plow$ \textit{undominated}. Any coherent lower probability (Definition~\ref{def:coherentlowerprob}) is dominated, but being dominated does not automatically imply coherence (\cf \citet[Section 2.4, Section 2.6.3]{walley1991statistical}).

 Write $D\coloneqq\{\omega^\infty: \exists A \subset \Omega: \limsupn r_n^{\omega^\infty}(A) >\liminf_{n \to\infty} r_n^{\omega^\infty}(A)\}$ for the event of divergence for at least one $A \subset \Omega$.
 \citet{kumar1985stationary} have shown the following (here presented with adapted notation\footnote{\citet{kumar1985stationary} use the letter $\mathcal{A}^\infty$ for what is $\infprodsigma$ in our notation, \ie the $\sigma$-algebra generated by all projection maps, or equivalently by the rectangles or cylinder sets.}).
\begin{proposition}[{\citep[Corollary 4.4]{kumar1985stationary}}]
    If $\Plow^\infty$ is monotone set function on $(\Omega^\infty,\infprodsigma)$ which is dominated, continuous from below on $\mathcal{A}^\infty$\footnote{Meaning if $\forall i \in \mathbb{N} : A_i^\infty \in \mathcal{A}^\infty$ and $A_i \uparrow A^\infty \in \sigma(\mathcal{A}^\infty)$, then $\Plow^\infty(A_i^\infty) \uparrow \Plow^\infty(A)$. Here, $\uparrow$ signifies ``non-decreasing'' and convergence of sets refers to the set-theoretic limit.} and stationary, then $\Plow^\infty(D)=0$. 
\end{proposition}
This result is intimately linked to the continuity condition. The reason for introducing this condition is that it is directly linked to the possibility of obtaining ``weak law''-type statements, which allow a bridge from finite observations to divergence or convergence. Indeed, the authors give examples very similar to (special cases of) our Definition~\ref{def:datamodelwithmnimprecision} of stationary coherent lower probabilities, hence dominated. In light of this result, however, such a model cannot be continuous from below. Intuitively, this prohibits any form of estimating cluster points from finite observations. As a conclusion, \citet{kumar1985stationary} suggest that set functions which satisfy the desired properties must be sought in the area of \textit{undominated} upper and lower probabilities, undominated implying \textit{incoherent} in the sense of \citet{walley1991statistical}.
The subsequent paper of \citet{papamarcou1991stationarity} provided an affirmative answer to the problem of finding a lower probability which is monotonely continuous along $\mathcal{A}^\infty$\footnote{Similar to continuity from below, see \citep{papamarcou1991stationarity}}, stationary, and has $\Plow(D)=1$ (a partial answer was provided in \citep{grize1987continuous}). In light of the above result, it is however undominated, hence incoherent. The final paper of \citet{sadrolhefazi1994finite} added to the result of \citet{papamarcou1991stationarity} by introducing an additional constraint, and proving that the problem is still feasible: the constraint that the lower probability should, on finite dimensional cylinder sets, behave approximately like an additive probability.

Philosophically, \citet{fine1988lower} subscribe to a propensity interpretation of the lower probabilities in these models. Due to the assumption of stationarity, the whole set of probabilities is considered to be ``at play'' for each ``timestep'' --- each ``timestep'' is characterized by an inherently imprecise propensity. To us, it remains unclear how to interpret such a propensity in the context of a data generation process. This clearly contrast this with the non-stationary, locally precise model in Section~\ref{sec:nonstationarymodels}, where each ``timestep`'' is linked to a single precise probability, and the imprecision manifests from a perspective of non-stationarity as aggregate (ir)regularity and local (ir)regularity. While we believe there might be applications for such more exotic data models like $\Plow_{\MtoNset}$ (\eg for certain physical phenomena such as flicker noise; see also Section~\ref{sec:applications}), we are mainly motivated by problems that are better described by non-stationarity. By going beyond stationarity, we can keep coherence and continuity, and still have models which typically yield divergence. Still there might be value in stationary models to describe a non-stationary process, since the stationary model might be more parsimonious, see \citep{persiau2021remarkable} for this argument. For example, non-stationary precise Markov chains can sometimes be efficiently modelled by stationary imprecise Markov chains \citep{de2009imprecise,t2019recursive}.

\subsection{A Comparison to Generalized Laws of Large Numbers and the Subjectivist Perspective}
\label{sec:llncomparisonsection}
 It is insightful to contextualize these data models in another vein of literature on generalized \textit{laws of large numbers} for various non-additive set functionals, going under the names of capacities or imprecise probabilities, and for the generalized expectation functionals called coherent lower (upper) previsions, coherent risk measures, or more recently, G-expectations \citep{dow1993laws,marinacci1999limit,epstein2003iid,maccheroni2005strong,peng2007law,de2008weak,rebille2009law,cozman2010concentration,chen2013strong,teran2014laws,hu2016general,peng2019nonlinear,zhang2020strong,zhang2024conditional}. We refer the reader to \citep{frohlich2022risk} for an account that explicates relations of capacities, coherent lower (upper) previsions and coherent risk measures.
We first describe the character of the mathematical results in these works, independent of interpretation. Then, since these authors primarily take a subjectivist perspective, it is also instructive to compare the interpretation of the mathematical objects, which may fulfill entirely different conceptual roles.

These results often take a certain form, which we informally describe to emphasize commonality over differences. Let $\Plow^\infty$ be a set functional which satisfies certain desirable properties, \eg is a convex capacity and continuous, and $\Rlow$,$\Rup$ be a corresponding pair of lower and upper generalized expectation\footnote{Often the Choquet integral is considered as the corresponding notion of generalized expectation. }. For concreteness, let $\Plow$ a coherent lower probability, obtained by applying the coherent lower prevision $\Rlow$ to indicator functions, \ie $\Plow(A) = \Rlow(\chi_A)$.
Let $X_i \sim X$ be a sequence of ``independent'' (in an appropriate \textit{imprecise} sense) gambles, then
\[
\Plow^\infty\left(\left\{\omega^\infty: \Rlow(X) \leq \liminf_{n\to\infty} \frac{1}{n} \sumin X_i(\omega^\infty) \leq \limsupn \frac{1}{n} \sumin X_i(\omega^\infty) \leq \Rup(X)\right\}\right)=1.
\]
Such a result is then called a \textit{strong law}. Similarly, a \textit{weak law} is roughly of the form
\[
\limn \Plow^n\left(\left\{\omega^\infty: \Rlow(X)-\varepsilon \leq \liminf_{n\to\infty} \frac{1}{n} \sumin X_i(\omega^\infty) \leq \limsupn \frac{1}{n} \sumin X_i(\omega^\infty) \leq \Rup(X)+\varepsilon\right\}\right)=1.
\]
For example, $X_i \coloneqq X \circ \pi_i$. Hence these statements are highly similar to what we obtained in Section~\ref{sec:comparisontofine}. As an example for an independence notion, consider Proposition~\ref{prop:modelonrectangles}, expressing a factorization on rectangles for $\Plow_{\MtoNset}$; in the literature, a variety of independence notions has been considered.\footnote{
The reader may consult the works cited above on generalized laws of large numbers to learn about the concrete independence concepts that they employ.}  Note however that such statements are silent on divergence or convergence within the interval (recall the discussion after Proposition~\ref{prop:wfweaklaw}), so that this kind of weak law is not in contradiction to the previously discussed impossibility of obtaining ``weak law''-type statements \textit{concerning apparent divergence} \citep{kumar1985stationary} for stationary, appropriately continuous and dominated lower probabilities.

While some authors emphasize the neutrality of their mathematical results \citep{de2008weak}, an often used interpretation in this literature is a subjectivist (generalized Bayesian) one, in the sense of \citet{gilboa1989maxmin} and \citet{walley1991statistical}, which can be viewed as a generalization of the approach by \citet{de2017theory}. Here, imprecision is taken to represent \textit{ambiguity}, sometimes called \textit{Knightian uncertainty}.
Ambiguity refers to a situation where a decision maker has no knowledge of objective probabilities and is unwilling to make precise, probabilistic assessments \citep{trautmann2015ambiguity,buhren2023ambiguity}.
We illustrate this will an \textit{Ellsberg example} \citep{ellsberg1961risk}, following the distinction of \citet{dow1993laws}.
Assume a decision maker who is \textit{uncertain} about the composition of urns, from which balls are drawn at random. It is known that the proportion of black balls is between $30\%$ and $40\%$, and the rest are red balls. Not intending to commit further, the decision maker might assign the lower $\Plow(B)=0.3$ and upper probability $\Pup(B)=0.4$ for drawing a black ball, and similarly $\Plow(R)=0.6$ and $\Pup(R)=0.7$. That is, the decision maker entertains a non-additive belief, which can be represented as a set $\mM$ of probability measures.
\citet{dow1993laws} distinguishes two long-run scenarios:
\begin{enumerate}[nolistsep,start=1,label=\textbf{U\arabic*.}, ref=U\arabic*]
    \item \label{item:rinftydatamodel} Balls are drawn independently (with replacement) from the same urn.
    \item \label{item:ourdatamodel} Balls are drawn independently from a sequence of urn, with possibly different compositions, but satisfying the above constraint on the proportions. There is no further evidence about the urn's compositions, implying that the urns are viewed as \textit{indistinguishable} by the decision maker \citep{epstein2003iid}.
\end{enumerate}
An appropriate subjective model for \ref{item:rinftydatamodel} is an extension to $\infprodsigma$ by a model proposed by \citet{walley1982towards}:
\[
\Plow^\infty_{\mM;\text{iid}}(A^\infty) \coloneqq \inf\{\operatorname{IT}((m,m,m,..))(A^\infty) : m \in \mM\}
\]
That is, this model is the envelope of all \iid product measures, but where any $m\in \mM$ is taken into account as a possibility.
In contrast, an appropriate subjective model for \ref{item:ourdatamodel} is given by our $\Plow_{\mM}^\infty$, which embodies a notion of independence, but leaves open the possibility that the compositions of the urns vary within $\mM$. Note that the above discussed strong laws, due to their independence assumptions, can be viewed as referring to situation~\ref{item:ourdatamodel} (see the Ellsberg example in \citep{chen2013strong}) rather than~\ref{item:rinftydatamodel}. The indistinguishability is transformed into a formal stationarity of the subjective model, even though the decision maker acknowledges that the true data generating process is non-stationary. Note that $\Plow^\infty_{\mM;\text{iid}}$ does not factorize in the same way as $\Plow_{\mM}$, thus it is not suitable for an imprecise independence concept (intuitively, according to $\Plow^\infty_{\mM;\text{iid}}$, if a probability measure $m \in \mM$ operates on the first trial, it also has to operate on all other trials).

This way of re-interpreting $\Plow_{\mM}$ from a subjectivist perspective makes clear the different conceptual role that it plays there. In situation \ref{item:ourdatamodel}, the actual data generating process is described by a non-stationary, locally precise data model of the form of Definition~\ref{def:preciseindependentdatamodel}. But since the decision maker is uncertain about the data model, they represent their uncertainty using a set of such data models, yielding $\Plow_{\mM}$. What is expressed is \textit{belief about the data}. Hence the models considered throughout Section~\ref{sec:stationarymodels} are suitable for different scenarios when looking through frequentist or subjectivist glasses.

Yet another way of interpreting the $\Plow_{\MtoNset}$ family of data models is this. Since it is constructed as the envelope of all \NSLP \xspace data models with aggregate (ir)regularity $\mathcal{N}$, the typicality $\Plow_{\MtoNset}(A^\infty)=1$ expresses that some event $A^\infty$ is typical \textit{with respect to any} of the respective \NSLP \xspace data models. This can then formally be used to justify that some estimation procedure (earning method) works for any such data model. In this sense, the typicality statement acts like a ``$\forall$'' quantifier, instead of aiming to describe the actual data generating process. This viewpoint however is in contrast to our attitude towards data models, which we interpret as describing the data generating process, untouched from any epistemic uncertainty about the data generating process itself. Compare this to how the standard \iid model is employed: in practice, we do not know the underlying probability $P$, but we use the \iid data model with respect to a single, but arbitrary $P$ as an abstract tool to justify our procedures, and we like to show that they work for any $P$.

\section{Applications and Estimability}
We have shown that, from a frequentist perspective, imprecision manifests in two distinct, yet interrelated, ways: as aggregate irregularity and as local regularity. Which of these manifestations we are interested in depends on context, essentially the decision problem at hand. Relatedly, these manifestations are linked to distinct ways of estimation (and thus learning problems). Before we discuss estimation, however, we believe it is insightful to first go into possible applications for these data models. The various, distinct such application contexts then shed light on which kind of estimation questions need to be asked and how they might be answered.

\subsection{Applications}
\label{sec:applications}
The non-stationary, locally precise (\NSLP) model keeps the familiar probabilistic concept, but departs from the \textit{identical} assumption. As such, we believe this model is an appropriate description of data-generating processes in many contexts. In light of the discussion in Section~\ref{sec:llncomparisonsection}, the models $\Plow_{\mM}$ and $\Plow_{\MtoNset}$ can describe subjective beliefs in such situations, for a decision maker who is uncertain about the data-generating process; the stationarity of the model would then describe this uncertainty.
Non-stationary application scenarios may include general data corruptions \citep{iacovissi2023general}, in particular dataset shift and learning with outliers, collaborative PAC learning, multi-source adaptation, fair machine learning and performativity. These problems share the common feature that they can be modelled by a set of distributions instead of a single one. Such multi-distribution settings have been receiving increasing attention in mainstream machine learning literature, see \eg \citep{haghtalab2022demand}.  For these problems, assume we work on a joint space $\Omega=\mathcal{X}\times \mathcal{Y}$, where $X: \Omega \to \mathcal{X}$ corresponds to the features of the learning problem, and $Y: \Omega \to \mathcal{Y}$ to the labels. We then consider a set of probability measures on $\Omega$, and seek to obtain a predictor $f(X)$ for $Y$.
\begin{enumerate}[nolistsep,start=1,label=\textbf{A\arabic*.},ref=A\arabic*]
    \item In dataset shift, we want to model the possiblity that the training and test distribution can differ: as machine learning systems are deployed in a dynamically unfolding environment, we can hardly expect stationarity. For safety reasons, it is therefore of great importance to robustify such systems against shifts in the underlying distribution.
    \item In the field of robust statistics \citep{huber1981robust}, neighborhood models have been widely studied to robustify statistical methods to outliers. 
    For example, as a frequentist version of the \textit{$\varepsilon$-contamination} model, assume that the data-generating process is described by an \NSLP \xspace data model parameterized by a sequence of probability measures $m^\infty=(P_0,..,P_0,..,P,..,P_0,..,P,..,P_0)$, where $P \in \Delta^k$ is arbitrary, and $P_0$ occurs with frequency at least $(1-\varepsilon)$ in the sequence in the sense that $\liminf_{n \to \infty} \chi_{\{j : P_j=P_0\}}(i) \geq 1-\varepsilon$. Intuitively, $P_0$ represents ``clean'' data and $P$ represents outliers (``corrupted'' data). A subjectivist decision maker who is uncertain about $P$ might then use the set $\mM = \{(1-\varepsilon)P_0 + \varepsilon P :  P \in \Delta^k\}$.
    \item In collaborative PAC learning \citep{blum2017collaborative} and federated learning \citep{mcmahan2017communication,zhang2021survey}, the goal is to find a predicator that performs well on multiple distributions, each corresponding to the data distribution of one client.
    \item In fair machine learning, the literature describes multiple competing goals \citep{barocas2023fairness}.  One among them is to obtain predictors which have equal predictive performance over multiple \textit{sensitive (ethically salient) subgroups}, determined for instance by gender, race \etc \citep{williamson2019fairness}. Each subgroup can be understood as corresponding to one probability measure on $\mathcal{X} \times \mathcal{Y}$, when the sensitive feature itself is not contained in $\Omega$.
    \item Under performativity \citep{perdomo2020performative}, the machine learning model itself supposedly acts on the data-generating process, yielding dataset shift in response to model deployment and therefore non-stationarity.
\end{enumerate}
The focus in these scenarios is often on the manifestation of local (ir)regularity (non-randomness, hidden heterogeneity), in the sense that the goal is often to obtain robust performance over the set of probabilities $\mM$. For instance, in the fair ML problem, hidden heterogeneity corresponds to the different ethically salient subgroups and the goal is not simply good performance on the whole population, but also on each of these subgroups. Sometimes, however, aggregate performance becomes the focus: for instance in outlier robustness, dataset shift or performativity we aim for overall good performance (on the whole aggregate) but acknowledge that a single, precise probability is insufficient as a data model.

The non-stationary, locally precise data model is couched in the language of precise probability. However, we view it also as inherently imprecise from the perspective of $\mM$ as a summary of the model. Is there any virtue to this view, or in different words, why would we want to work with imprecise probabilities in such settings? By this, we mean for example introducing predictors which provide imprecise probabilistic forecasts as opposed to precise ones. We think there are legitimate reasons for doing so. \citet{walley1982towards} have offered the following intuition for the usefulness of imprecise modelling:
\begin{quote}
    ``Just as `randomness' (chance) is introduced in additive probability models to account for poorly understood (`accidental') variation in outcomes, so `indeterminacy' might be introduced in upper and lower probability models to account for poorly understood variations in chance behaviour.''
\end{quote}
In this way, imprecision (here called ``indeterminacy'') relates to a second-order modelling level, distinct from first-order randomness. In the standard \iid picture, the goal would be to approximate the ideal Bayes predictor $f(X)=\mathbb{E}_P[Y|X]$ with respect to the precise underlying probability $P$. Intuitively, all information in $X$ should be used by the predictor. However, the above problems make plausible that this is not always the case: in collaborative PAC learning, we know which distribution corresponds to which agent, but we might not want our predictor $f$ to depend on this. Similarly, in a fairness context, it might be undesirable or indeed prohibited to include the sensitive feature in $\mathcal{X}$. In this way, the above problems share a \textit{hidden heterogeneity} which cannot or should not be captured in $\mathcal{X}$ (which would reduce the problem to a single, precise distribution), but with robustness desiderata with respect to the hidden heterogeneity. We believe that the main use of imprecision for learning under an \NSLP\xspace data model is to account for ``poorly understood variations in chance behaviour'' as claimed above, but still understood \textit{to some extent}. That is, there must be some handle available to get to this variation (see the discussion on estimability below).

Besides its subjectivist interpretation, what are possible applications for $\Plow_{\mM}$ and $\Plow_{\MtoNset}$ as data models? We have already noted that they might serve as more parsimonious stationary models for precise, non-stationary phenomena. For example, when there is no understanding \textit{at all} of such variation in chance, in contrast to above. However their actual character would prove to be relevant when we assume indeterminacy in chance \textit{even under stationarity}. Applications under this philosophy remain more speculative; \citet{fine1988lower} have given flicker noise as an example for a process which is hypothesized to be stationary but exhibits aggregate irregularity. 
We believe the prime area of application might be \textit{heavy-tailed phenomena}, which are recently receiving increasing attention \citep{taleb2007black,resnick2007heavy,nair2022fundamentals} --- already \citet{fierens2001towards}, who focus on aggregate irregularity, suggested ``[..] that these models will be
appropriate when dealing with such long sequences as might arise from heavy-tailed waiting times.''. For heavy-tailed phenomena, an often-used model is that of a stationary Pareto distribution with parameter $\alpha$, which produces relative frequencies that converge on much slower rates than for light-tailed distributions for $\alpha>1$; for $\alpha<1$, relative frequencies diverge. However, \citet{frohlich2024strictly} have argued that divergence of relative frequencies in the infinite limit is merely an idealization of unstable relative frequencies (``apparent divergence'' as in Section~\ref{sec:comparisontofine}) for finite $n$, in the same way that convergence idealizes apparent stability; in this way, a data model which typically yields divergence may also be appropriate for $\alpha>1$. We leave open the question of whether (or which) such heavy-tailed phenomena are better understood and modelled as non-stationary, locally precise or stationary, locally imprecise. In fact, recent results in imprecise game-theoretic probability support the claim that there is a fundamental equivalence relationship between non-stationary, locally precise and stationary, locally imprecise models \citep{persiau2021remarkable,cooman2022randomness}. From a practical point of view, a relevant area of applications would be machine learning and statistics with heavy-tailed outcomes. Such use cases arise for instance in (re)insurance uses cases, where claims (losses) can exhibit heavy-tails \citep{ibragimov2015heavy,peng2017inference,PUNZO201895}. As another interesting research avenue, recent publications in the field of climate science \citep{lovejoy2015voyage,franzke2020structure} uncover aggregate irregularity in climate time series. We take this to suggest the development of a foundational account which links together and explicates the relationships between the concepts of power laws, heavy-tails, (non)-stationarity and aggregate irregularity.

 Given the manifold nature of these application scenarios are, we shy away from constructing more concrete, constrained data models in the present paper. Rather, our aim is to offer a firm foundation on which such subsequent work can take place. We do however want to make some general remarks on estimation in such data models, which then needs to be tailored to the specific scenarios.

\subsection{Estimation of Aggregate (Ir)regularity}
Assume a fixed sequence $\omega^\infty=(\omega_1,\omega_2,..)$. Can we estimate its cluster points of relative frequencies, \ie $\CP(n \mapsto r_n^{\omega^\infty})$? In fact, there seems hardly a reason for attempting to do this. This set could potentially be uncountable, which would make for a daunting task. Also observe that knowing $\Plow(\{\omega^j\})$, $\Pup(\{\omega^j\})$ for each elementary event $\omega^j$ is not sufficient for knowing whole set of cluster points, since a coherent lower (upper) probability is not in a one-to-one correspondence to a set of measures.\footnote{Recall that a coherent lower (upper) prevision is in a one-to-one correspondence to a closed convex set of (finitely additive) probabilities.}

Further, in practice we do not observe elementary events, but rather the values that specific gambles take. In a decision problem, we are rather interested in estimating quantities such as
\[
\liminf_{n \to \infty} \frac{1}{n} \sumin \ell(\omega_i) \quad \text{ and } \limsupn \frac{1}{n} \sumin \ell(\omega_i)
\]
for a loss gamble $\ell: \Omega \to \mR$ (\eg a loss function partially applied to an action). This is a much simpler problem, and indeed can be solved to some extend. Note that due to connectedness, we have 
\[
\CP\left(n \mapsto  \frac{1}{n} \sumin \ell(\omega_i)\right) = \left[\liminf_{n \to \infty} \frac{1}{n} \sumin \ell(\omega_i),\limsupn \frac{1}{n} \sumin \ell(\omega_i)\right]
\]
Thus, for gamble averages, an estimator for the $\liminf$ would due to conjugacy suffice to estimate the whole interval of cluster points. First, we observe that this problem is in principle solvable for a fixed sequence.
\begin{proposition}
    For any $\omega^\infty$, there exists functions $g_n$ of the form $g_n(\ell(\omega_1),..,\frac{1}{n} \sumin \ell(\omega_i))$ so that 
    \[
    \limn g_n(..) = \liminf_{n \to \infty} \frac{1}{n} \sumin \ell(\omega_i).
    \]
\end{proposition}
\begin{proof}
    $\mR$ with the Euclidean topology is sequentially compact, and thus each cluster point of the sequence of gamble averages corresponds to the limit point of a convergent subsequence. Hence there exists a selection rule $S: \mathbb{N} \to \{0,1\}$ which extracts this subsequence .
\end{proof}
Not only does this ``estimator'' depend on the concrete sequence $\omega^\infty$, it is completely useless, since the selection rule is unknown in practice. Using ``all'' selection rules, on the other hand, would not meaningfully work since in the finite case this would imply the ``vacuous'' estimate $\{\min_{1 \leq j \leq n} \frac{1}{j} \sum_{i=1}^j \ell(\omega_i), \max_{1 \leq j \leq n} \frac{1}{j} \sum_{i=1}^j \ell(\omega_i)]$.
We therefore desire an estimator which works on all sequences $\omega^\infty$.
\citet{walley1982towards} have proposed the following estimator:\footnote{\citet{walley1982towards} proposed the estimator for relative frequencies of events, but the generalization to gamble averages is obvious. Originally, they demanded that $\limn \kfunc(n)/n = 0$, but the proof of this statement in fact works if $\kfunc(n)\leq n$.}
\begin{equation}
    \label{eq:wfestimator}
    \hat{\underline{\operatorname{avg}}}_{\ell,n} = \min\left\{\frac{1}{j} \sum_{i=1}^j \ell(\omega_i) : \kfunc(n) \leq j \leq n\right\},
\end{equation}

where the function $\kfunc : \mathbb{N} \to \mathbb{N}$ is such that $\limn \kfunc(n)=\infty$ and $\limn \kfunc(n) \leq n$. For example, $\kfunc(n)=\sqrt{n}$. Similarly, $\hat{\overline{\operatorname{avg}}}_{\ell,n} $ is defined using $\max$ instead of $\min$. 
\begin{proposition}
\label{prop:weakestimation}
    This estimator succeeds in the sense that:
\[
\liminf_{n \to \infty} \avgest = \liminf_{n \to \infty} \frac{1}{n} \sumin \ell(\omega_i) 
= \min\left\{\mathbb{E}_p[\ell] : p \in \CP(n \mapsto r_n^{\omega^\infty})\right\}.
\]
\end{proposition}
\begin{proof}
    The second equality of the statement is due to Proposition~\ref{prop:cpsofrelativefrequenciesandgambles}.
    The first equality follows straightforwardly by adapting the proof of Theorem 4.2 by \citet{walley1982towards}. For completeness: for fixed $m \in \mathbb{N}$, $\inf\{\avgest : n \geq m\} = \inf\{\frac{1}{n} \sumin \ell(\omega_i)  : n \geq \min_{n \geq m} \kfunc(n)\} \eqqcolon \RomanNumeralCaps{1}(m)$. Since $\kfunc(n) \to \infty$  we then get that $\lim_{m \to \infty} \RomanNumeralCaps{1}(m) = \lim_{m \to \infty} \inf\{\frac{1}{n} \sumin \ell(\omega_i) : n \geq m\} =  \liminf_{n \to \infty} \frac{1}{n} \sumin \ell(\omega_i)$.
\end{proof}
An analogous statement holds for the conjugate upper quantity. In this sense, the estimation succeeds asymptotically, but in a very weak sense only: we only get that the $\liminf$ of the estimator gives the desired estimand. This result can hardly be considered satisfying; also note that $\kfunc(n)=n$ in the above simply gives the familiar estimator $\avgest =  \frac{1}{n} \sumin \ell(\omega_i)$, which has the same guarantee due to Proposition~\ref{prop:weakestimation}. The following suggests that \textit{no choice} of $\kfunc(n)$ will work in general in the stronger sense that the limit of the estimator would coincide with the estimand.
\begin{proposition}
\label{prop:veryslowsequence}
    For any choice of $\kfunc(n)$ which satisfies $\limn \kfunc(n) = \infty$ and $\kfunc(n)\leq n$ and set of probability measures $\mathcal{M} \subseteq \Delta^k$, $|\mathcal{\mM}|>1$, and any $\mathcal{N} \in \CM$, there exists a sequence of measures $m^\infty \in \mM^\infty$ so that
    \begin{equation}
    \label{eq:wfcps}
    \CP\left(n \mapsto \min\left(\frac{1}{j} \sum_{i=1}^j m_i : \kfunc(n) \leq j \leq n\right) \right) = \mathcal{N}.
    \end{equation}
\end{proposition}
\begin{proof}
    Note that it is clear that $\CP\left(n \mapsto \min\left(\frac{1}{j} \sum_{i=1}^j m_i : \kfunc(n) \leq j \leq n\right) \right) \subseteq \CP\left(n \mapsto \frac{1}{n} \sumin m_i \right)$. We now want to show that the sequence $m^\infty$ can be constructed so that the set on the left side is not a singleton, \ie the limit does not exist. Proposition~\ref{prop:forwarddir}, shouldering the main work of the proof of Theorem~\ref{theorem:maintheoremffivanenko}, has shown the existence of some $m^\infty$ with the property that $\CP\left(n \mapsto \frac{1}{n} \sumin m_i \right)=\mathcal{N}$ for any $\mathcal{N} \in \CM$. We now show that the proof of Proposition~\ref{prop:forwarddir} can be slightly modified (or more precisely, can be concretly instantiated) so that \eqref{eq:wfcps} holds in addition. Intuitively, the average $\frac{1}{n} \sumin m_i$ has to change extremely slow relative to $\kfunc(n)$, so that $\{\frac{1}{j} \sum_{i=1}^j m_i : \kfunc(n) \leq j \leq n \}$ is approximately constant. The proof of Proposition~\ref{prop:forwarddir} contains two relevant ``gaps'', where it is only asserted that there exists a finite $\kfunc \in \mathbb{N}$ will make a condition hold, but no concrete $\kfunc$ is chosen:
    \begin{enumerate}
        \item \eqref{eq:first_implementation_gap}: $\exists k \in \mathbb{N}: r(m \ccm b^k) \in B_{\kappa}(\tilde{c}_2)$,
        \item \eqref{eq:second_implementation_gap}: $\exists k \in \mathbb{N}: r(m \ccm b^k) \in B_\kappa(q).$
    \end{enumerate}
    The proof therefore needs to be instantiated with a concrete choice of $k$ in these statements. Note that if either condition holds for some $k$, it also holds for all $l \geq k$ due to Lemma~\ref{lemma:relfreqsconvex}. It is guaranteed that there exists always a choice of $k$ so that the proof yields $\CP\left(n \mapsto \frac{1}{n} \sumin m_i \right)=\mathcal{N}$. But further, if we choose $k$ large enough, we can in addition guarantee that \eqref{eq:wfcps} will hold for the resulting sequence $m^\infty$: we need to find some $k^*$ so that from $j=\kfunc(k^*)..k^*$, $r(m \ccm b^{j}) \in B_\kappa(\tilde{c}_2)$ (analogously for \eqref{eq:second_implementation_gap}), which is clearly possible. This implies that $
     \min\left(\frac{1}{j} \sum_{i=1}^j m_i : \kfunc(k^*) \leq j \leq k^*\right) \in B_\kappa(\tilde{c}_2)$. The logic of the proof of Proposition~\ref{prop:forwarddir} then implies \eqref{eq:wfcps}.
\end{proof}

We conclude by applying Proposition~\ref{prop:samecpastheoretical} that with this estimator, not even the cluster points of relative frequencies can be estimated, dispensing with the hope to estimate gamble averages. We now write $\hat{\underline{\operatorname{avg}}}_{\ell,n}^{\omega^\infty}$ to highlight the dependence on the sequence $\omega^\infty$.
\begin{corollary}
\label{cor:wfestimatorfails}
Let $\emptyset \subsetneq A \subsetneq \Omega$. 
    For any choice of $\kfunc(n)$ with $\limn \kfunc(n) = \infty$ and $\kfunc(n)\leq n$, there exists some sequence $\omega^\infty$ so that 
    $\lim_{n \to \infty} \avgesta^{\omega^\infty}$ does not exist.
\end{corollary}
\begin{proof}

Choose probability measures $\mM \coloneqq \{m_1,m_2\}$ so that $m_1(A) \neq m_2(A)$. Proposition~\ref{prop:samecpastheoretical} yields that the sequence $m^\infty \in \mM^\infty$ gives $\operatorname{IT}(m^\infty)(\{\omega^\infty: \lim_{n \to \infty} \avgesta^{\omega^\infty} \text{ does not exist}\})=1$, therefore one such sequence exists.
\end{proof}
Until now we considered arbitrary sequences $\omega^\infty$.
To contextualize estimation in our data models, we emphasize that estimation is inextricably tied to typicality: only for typical sequences can we hope to recover the underlying parameterization; for example, even a familiar \iid model can generate sequences with diverging relative frequencies, but only \textit{atypically}. The previous result suggests that this estimator will not \textit{typically} work in our data models. First, however, we look at a positive result. \citet{walley1982towards} have shown that, if $\kfunc(n)$ satisfies a specific growth behaviour, estimation works at least in the sense of asymptotic favorability. Recall that $A_1, A_2, ..$ are asymptotically favored under $\Plow$ when $\limn \Plow((A_n)^c)  / \Plow(A_n) = 0$, implying $\limn \Pup(A_n)=1$. The result is stated in terms of relative frequencies instead of gamble averages.
\begin{proposition}{{\citep[Theorem 5.3]{walley1982towards}}}
    Assume $\kfunc(n)/n \to 0$ and $\kfunc(n)/\log(n) \to \infty$. Let $\Pup(A)<1$ or $\Plow(A)=1$. Then
    \[
    \forall \varepsilon>0:  \{\omega^\infty: |\avgesta^{\omega^\infty} - \Plow_{\mM}(A)|<\varepsilon\} \text{ asymptotically favored under } \Plow_{\mM}^\infty.
    \]
\end{proposition}
Since they agree on $\mathcal{A}^\infty$, the same holds for $\Plow_{\MtoNset}^\infty$. On the other hand, it is easily seen that if $0< \varepsilon<\Pup_{\mM}(A)-\Plow_{\mM}(A)$ then $\Plow_{\mM}^\infty(\{\omega^\infty: |\avgesta^{\omega^\infty}  - \Plow_{\mM}(A)|<\varepsilon\}) = 0$ and therefore a ``weak law''-type statement for estimation fails. However, from Proposition~\ref{prop:veryslowsequence} we obtain:
\begin{proposition} Let $\emptyset \neq A \subsetneq \Omega$ an imprecise event, that is, $\Pup(A)>\Plow(A)$. For any choice of $\kfunc(n)$ with $\limn \kfunc(n) = \infty$ and $\kfunc(n)\leq n$, the Walley-Fine estimator fails in the sense that
    \[
    \Plow_{\MtoNset}\left(\left\{\omega^\infty : \limn \avgesta^{\omega^\infty} \text{ exists } \right\}\right) = 0,
    \]
    and similarly for $\Plow_{\mM}^\infty$, since $\Plow_{\mM}^\infty \leq \Plow_{\MtoNset}^\infty$.
\end{proposition}
\begin{proof}
    Observe that the sequence $m^\infty=(m_1,m_2,..)$ of Proposition~\ref{prop:veryslowsequence} lies in $\MtoNset$.
\end{proof}
This result seems unsettling, and motivates us to conjecture the following.
\begin{conjecture}
    There is no ``reasonable'' (without access to an ``oracle'') estimator so that
    \[
    \Plow_{\MtoNset}\left(\left\{\text{ estimation of cluster points of relative frequencies succeeds } \right\}\right) = 1.
    \]
\end{conjecture}

On the other hand, for \textit{some} sequences of probabilities, the estimator \textit{will} work typically.
\citet[p.\@\xspace 13]{kumar1985stationary}, with a similar motivation, have therefore reversed the logic and proposed to construct the model exactly in a way such that the estimator works in a ``strong law'' sense (by collecting all sequences of probabilities for which the estimator with a specified $\kfunc$ succeeds almost surely and forming the envelope). We believe this is indeed the way forward: in our view, $\Plow_{\mM}$ and $\Plow_{\MtoNset}$ only represent the starting point for developing data models which are finely tailored to specific assumptions about rates of convergence (divergence).

On a broader level, do we even \textit{want} to estimate aggregate irregularity? Our position is that when considering estimation we cannot ignore the relevant decision-making context, as it determines which questions are meaningful to be asked. We believe that an explicit estimation of aggregate irregularity is only called for when the decision problem is with respect to aggregate criteria \textit{and} context demands stationary (constant) forecasts. Recall the criticsm by \citet{fierens2003towards} discussed in Section~\ref{sec:fierensfinecomparison}: if a skeptic can probe for local regularity, a forecast based on aggregate irregularity may perform poorly. Intutively, the problem of the estimator \eqref{eq:wfestimator} is that relative frequencies are always computed from the very beginning; therefore the estimator adapts extremely slowly to changing \textit{local} relative frequencies in the data sequence, \eg as computed by averages over moving windows, whereas a non-stationary forecast could quickly adapt to them. This shows that if evaluation is based on local (ir)regularity, estimation should also focus on this manifestation.
On the other hand, when the evaluation of our forecasts is based on the whole aggregate \textit{and} there is a demand for stationary (constant) forecasts, estimating aggregate (ir)regularity becomes meaningful. For instance, in the example in Section~\ref{sec:applications} of machine learning with heavy-tailed insurance claims, our evaluation is on the aggregate, and we might have no reason to suspect a time evolution and thus may want a stationary model. Finally, when evaluation is based on aggregate criteria, but there is no demand for stationary (constant) forecasts, a wider space opens up: an intruiging recent work by \cite{zhao2021right}, which exhibits remarkable parallels to imprecise probabilities, has provided an online learning method that, based on non-stationary forecasts, satisfies a strong performance guarantee even under aggregate irregularity. In other words, their scheme passes a test in terms of aggregate-level performance even under aggregate irregularity without directly aiming to estimate aggregate irregularity; likewise for other online learning methods \citep{vovk2005defensive}.

To summarize, what is important is arguably not to estimate aggregate irregularity per se, but to have learning methods which have performance guarantees even under aggregate irregularity, however this is achieved.

\subsection{Estimation of Local (Ir)regularity}
\label{sec:estimationoflocal}
First, assume a non-stationary, locally precise model parameterized by the sequence $m^\infty=(m_1,m_2..)$ so that $\mM\coloneqq \{m_1,m_2,..\}$. In problem contexts where the focus is on guarding against hidden heterogeneity, we may wish to estimate $\mM$. For this, assume we have access to a fixed, finite set of selection rules $\mathcal{S}$. For a selection rule $S \in \mathcal{S}$, define its induced subsequence as $m^\infty_S$ as the subsequence of $m^\infty$ extracted by the condition that $S(i)=1$.
The theoretical mean associated to a selection rule is: 
\[
\mu_{S,n}(\{\omega^j\}) \coloneqq \frac{\frac{1}{n} \sum_{i=1}^n m_i(\{\omega^j\}) S(i)}{\frac{1}{n} \sumin S(i)}, \quad \omega^j \in \Omega,
\]
(well-defined beginning with some $n_0 \in \mathbb{N}$ for which $S(n_0)=1$) which extends to a probability measure $\mu_{S,n}$ on $\Delta^k$ in the obvious way by additivity.
On the other hand, for a fixed sequence $\omega^\infty$, the empirical relative frequencies induced by a selection rule are as follows:
\[
r_{S,n}^{\omega^\infty}(\{\omega^j\}) \coloneqq \frac{\frac{1}{n} \sumin \chi_{\{\omega^j\}}(\omega_i) S(i)}{\frac{1}{n} \sumin S(i)},\quad \omega^j \in \Omega,
\]
(again, well-defined for all $n \geq n_0 \in \mathbb{N}$) 
which similar extends to a probability measure on $\Delta^k$.
\citet[Theorem 1]{fierens2009frequentist} have shown the following.
\begin{proposition}
    Let $\operatorname{D}(p,q) \coloneqq \max_{\omega \in \Omega} |p(\{\omega\}) - q(\{\omega\})|$. Let $0 < m\leq n$. Then:
    \[
    \lambda\left\{\omega^\infty \coloneqq G(z) : \max_{S \in \mathcal{S}}\left(D(\mu_{S,n},r_{S,n}^{\omega^\infty}) : \frac{1}{n} \sumin S(i) \geq m \right) \geq \varepsilon \right\} \leq 2 k |\mathcal{S}| \exp\left(-\varepsilon^2 m^2 / 2n\right)
    \]
\end{proposition}
Intuitively, the familiar weak law holds for each selected subsequence, and we can therefore avoid extracting ``arbitrary patterns'', but we must consider only those selection rules which have selected sufficiently many indices until $n$. For example, if $m^\infty$ truly was an \iid model, \ie $m^\infty=(m_0,m_0,..)$, the above statement reduces to a weak law, asserting that relative frequencies converge to $m_0$ along all selection rules. On the other hand, if there is underlying imprecision in $m^\infty$, manifested as non-stationarity, then selection rules might converge to different limits, or indeed they might not converge at all; Theorem~\ref{theorem:maintheoremffivanenko} implies that $\CP\left(n \mapsto \mu_{S,n}\right) \in \mathcal{C}(\mM)$.

A set of selection rules $\mathcal{S}$ therefore provides a natural estimate for $\mM$:
\[
\hat{M}_n^{\omega^\infty} \coloneqq \bigcup_{S \in \mathcal{S}} r_{S,n}^{\omega^\infty}
\]
For the purpose of computing lower (upper) expectations, the closed convex hull might be used equivalently. 
The implication is that for large enough $n$, we do not (approximately) over-estimate the set $\mM$, that is, our estimate is (approximately) contained in $\mM$. On the other hand, we have no guarantee that we will recover \textit{all of $\mM$}. Access to selection rules of the form of Proposition~\ref{prop:uncoverhiddenheterogeneity} would at least guarantee that we asymptotically recover all $m \in \CP(i \mapsto m_i)$. 
But assuming access to such a set seems highly impractical. Motivated by the goal of recovering the whole $\mM$, \citet{fierens2009frequentist} impose computability requirements, to argue that they can recover $\mM$ ``universally''.

The situation can be clarified by taking a step back and reflecting about how such data models may actually be used in practice. In practice, we do not have access to an oracle that gives us the ``true'' data model (we put aside the question of the ontological nature of such data models) but rather we start with an approach to data handling/processing. Concretely, we may construct a set of selection rules, guided by the application context, \ie starting with assumptions about how the hidden heterogeneity might be structured. Then we invoke a data model to argue, or justify, the sensibility of our approach, by showing that it has desirable properties under this data model. This line of reasoning sheds light also on the familiar \iid model: it is not because real world data is drawn \iid that the \iid model is so popular; rather, because it is customary to evaluate based on average performance and ignore the heterogeneity of subpopulations that the \iid model serves to justify such an approach. For a more general example, if in some context the relevant properties of our decision problems depend only on the aggregate (ir)regularity, there is no difference for whether we would use $\Plow_{\MtoNset}$ as a data model or an \NSLP \xspace data model with aggregate (ir)regularity $\mathcal{N}$. The choice of data model (crucially, its typicality notion) should reflect our data processing, in particular with respect to the relevant characteristics of the decision problem at hand.

On the other hand, there are of course limits to the ``self-fulfilling'' prophecy of assuming a data model. The data we have actually at hand might be \textit{atypical} with respect to our assumed data model. For example, we may have assumed an \iid model, but investigating with selection rules might uncover local regularity in the data, atypical for an \iid model.\footnote{Here we are using typicality in an approximate, non-asymptotic sense. For small, finite $n$, even an \iid model could yield local regularity, but we expect it to disappear for large enough $n$.} This should motivate us to then reject the \iid model and to look for a data model which is consistent with the data at hand --- this is also the logic of a frequentist hypothesis test.

For $\Plow_{\mM}$ we note that if $|\mM|>1$, the model is vacuous about the existence of local (ir)regularity, since the corresponding set of measures also contains \iid product measures. It seems hardly feasible to obtain interesting statements of the local (ir)regularity associated with $\Plow_{\MtoNset}$. If $\mathcal{N}$ is not a singleton, no \iid product measure can be in the corresponding set of measures, implying that each such product measure has \textit{some} local regularity. However, since the set of measures is only constrained by their divergent behaviour, it appears impossible that they ``agree'' on the structure of the local regularities --- otherwise the model could not be stationary.

\section{Discussion}
\citet{fine1970apparent} had already noted about aggregate regularity:
\begin{quote}
    \textit{Apparent
convergence of the relative frequency occurs because of, and
not in spite of, the high irregularity
(randomness)
of the
data sequence.} [emphasis in original]
\end{quote}
Our main theorem implies that we only partially agree with this finding: indeed, local regularity (non-randomness) is a prerequisite for aggregate irregularity (divergence). But the existence of local regularity as a hidden heterogeneity need not imply aggregate irregularity.
While we believe we have thus fully characterized the relation between these two manifestations of imprecision, and explicated how they manifest under different data models, many new questions have surfaced in this process.
One line of future research would develop more constrained data models with imprecision, tailored also to specific application contexts and provide methods for estimation. For instance, can we craft a notion of typicality directly in terms of rate of convergence (divergence)? In light of Section~\ref{sec:comparisontofine}, such models would necessarily be non-stationary, if they are to be coherent.

Often, imprecision is taken to be synonymous with subjective ambiguity, where the idea is that it reflects merely epistemic uncertainty. We have departed from this picture by introducing data models which cannot be described by a single probability, irrespective of our epistemic considerations. This suggests that there is a role for imprecise forecasts in this picture, where a predictor would be of the type $f(X=x) \in 2^{\Delta^k}$, outputting for each $X=x$ a set of probabilities. Under an imprecise data model like \NSLP \xspace or \SLI \xspace parameterized by $\mM$, it is plausible that forecasting the set $\mM$ itself is a ``reasonable'' stationary forecast; but how to evaluate IP forecasts? A rigorous development, in our opinion, would introduce and study generalizations of proper scoring rules and calibration to the imprecise case. Impossibility results in the literature \citep{seidenfeld2012forecasting,mayo2016scoring,schoenfield2017accuracy}, directed against IP scoring rules, \textit{prima facie} threaten this endeavour. We do not share the pessimism --- see also \citep{konek2019ip} for a discussion --- and highlight that these results are concerned with \textit{strict propriety}. Intuitively, this is the demand that there cannot, to any imprecise forecast, exist a precise forecast which is at least as accurate as measured by the scoring rule. In line with \citet{mayo2016scoring} we think that the strictness is an overly stringent desideratum. Preliminary findings which we have obtained indeed indicate the feasibility of developing weakly proper IP scoring rules and IP calibration. In our view, this requires relativizing these concepts to the assumed data model; in previous works on IP scoring rules which follow a subjectivist approach, data models have not appeared. The present paper provides the foundation for proceeding further in this direction.

Another follow-up question is how to represent epistemic uncertainty, when we use imprecision to represent aleatoric indeterminacy. In an estimation process like that of Section~\ref{sec:estimationoflocal} the estimate, which is a set of probabilities, is itself subject to epistemic uncertainty, which calls for a principled, disentangled representation of the distinct sources of uncertainty.

Finally, in the present paper we have not rigorously linked and compared these frequentist data models to recent work on imprecise game-theoretic probability \citep{persiau2021remarkable,persiau2022on,cooman2022randomness}. These authors also take inspiration from the works of \citet{walley1982towards} and subsequent papers by Fine and collaborators, but operate in the language of game-theoretic probability, rendering a translation of results non-trivial.

\section{Acknowledgements}
This work was funded by the Deutsche Forschungsgemeinschaft (DFG, German Research Foundation)
under Germany’s Excellence Strategy — EXC number 2064/1 — Project number 390727645. The authors
thank the International Max Planck Research School for Intelligent Systems (IMPRS-IS) for supporting
Christian Fröhlich. Many thanks to Rabanus Derr and Marvin Pförtner for helpful discussions.

\bibliography{twofacesofimprecision.bib}

\begin{thebibliography}{89}
\providecommand{\natexlab}[1]{#1}
\providecommand{\url}[1]{\texttt{#1}}
\expandafter\ifx\csname urlstyle\endcsname\relax
  \providecommand{\doi}[1]{doi: #1}\else
  \providecommand{\doi}{doi: \begingroup \urlstyle{rm}\Url}\fi

\bibitem[Artzner et~al.(1999)Artzner, Delbaen, Eber, and Heath]{artzner1999coherent}
Philippe Artzner, Freddy Delbaen, Jean-Marc Eber, and David Heath.
\newblock Coherent measures of risk.
\newblock \emph{Mathematical Finance}, 9\penalty0 (3):\penalty0 203--228, 1999.

\bibitem[Augustin(2022)]{augustin2022statistics}
Thomas Augustin.
\newblock Statistics with imprecise probabilities—a short survey.
\newblock In \emph{Uncertainty in Engineering: Introduction to Methods and Applications}, pp.\  67--80. 2022.

\bibitem[Augustin et~al.(2014)Augustin, Coolen, De~Cooman, and Troffaes]{augustin2014introduction}
Thomas Augustin, Frank~P.A. Coolen, Gert De~Cooman, and Matthias~C.M. Troffaes.
\newblock \emph{Introduction to imprecise probabilities}.
\newblock John Wiley \& Sons, 2014.

\bibitem[Barocas et~al.(2023)Barocas, Hardt, and Narayanan]{barocas2023fairness}
Solon Barocas, Moritz Hardt, and Arvind Narayanan.
\newblock \emph{Fairness and machine learning: Limitations and opportunities}.
\newblock MIT Press, 2023.

\bibitem[Blum et~al.(2017)Blum, Haghtalab, Procaccia, and Qiao]{blum2017collaborative}
Avrim Blum, Nika Haghtalab, Ariel~D. Procaccia, and Mingda Qiao.
\newblock Collaborative {PAC} learning.
\newblock In \emph{Advances in Neural Information Processing Systems}, volume~30. 2017.

\bibitem[Bradley \& Steele(2014)Bradley and Steele]{bradley2014should}
Seamus Bradley and Katie Steele.
\newblock Should subjective probabilities be sharp?
\newblock \emph{Episteme}, 11\penalty0 (3):\penalty0 277--289, 2014.

\bibitem[B{\"u}hren et~al.(2023)B{\"u}hren, Meier, and Ple{\ss}ner]{buhren2023ambiguity}
Christoph B{\"u}hren, Fabian Meier, and Marco Ple{\ss}ner.
\newblock Ambiguity aversion: bibliometric analysis and literature review of the last 60 years.
\newblock \emph{Management Review Quarterly}, 73\penalty0 (2):\penalty0 495--525, 2023.

\bibitem[Chen et~al.(2013)Chen, Wu, and Li]{chen2013strong}
Zengjing Chen, Panyu Wu, and Baoming Li.
\newblock A strong law of large numbers for non-additive probabilities.
\newblock \emph{International Journal of Approximate Reasoning}, 54\penalty0 (3):\penalty0 365--377, 2013.

\bibitem[Cozman(2010)]{cozman2010concentration}
Fabio Cozman.
\newblock Concentration inequalities and laws of large numbers under epistemic and regular irrelevance.
\newblock \emph{International Journal of Approximate Reasoning}, 51\penalty0 (9):\penalty0 1069--1084, 2010.

\bibitem[Cozman \& Chrisman(1997)Cozman and Chrisman]{cozmanconvex}
Fabio Cozman and Lonnie Chrisman.
\newblock Learning convex sets of probability from data.
\newblock Technical report, Carnegie Mellon University, 1997.
\newblock CMU-RI-TR 97-25.

\bibitem[De~Cooman \& De~Bock(2022)De~Cooman and De~Bock]{cooman2022randomness}
Gert De~Cooman and Jasper De~Bock.
\newblock Randomness is inherently imprecise.
\newblock \emph{International Journal of Approximate Reasoning}, 141:\penalty0 28--68, 2022.

\bibitem[De~Cooman \& Miranda(2008)De~Cooman and Miranda]{de2008weak}
Gert De~Cooman and Enrique Miranda.
\newblock Weak and strong laws of large numbers for coherent lower previsions.
\newblock \emph{Journal of Statistical Planning and Inference}, 138\penalty0 (8):\penalty0 2409--2432, 2008.

\bibitem[De~Cooman et~al.(2009)De~Cooman, Hermans, and Quaeghebeur]{de2009imprecise}
Gert De~Cooman, Filip Hermans, and Erik Quaeghebeur.
\newblock Imprecise markov chains and their limit behavior.
\newblock \emph{Probability in the Engineering and Informational Sciences}, 23\penalty0 (4):\penalty0 597--635, 2009.

\bibitem[de~Finetti(1974/2017)]{de2017theory}
Bruno de~Finetti.
\newblock \emph{Theory of probability: A critical introductory treatment}.
\newblock John Wiley \& Sons, 1974/2017.

\bibitem[Dow \& Werlang(1993)Dow and Werlang]{dow1993laws}
James Dow and S{\'e}rgio Ribeiro da~Costa Werlang.
\newblock Laws of large numbers for non-additive probabilities.
\newblock 1993.
\newblock URL \url{https://repositorio.fgv.br/server/api/core/bitstreams/5e83654c-e8a5-4ad7-b6e5-6ba93ee58a6f/content}.
\newblock Accessed: 2024-04-03.

\bibitem[Eagle(2021)]{sep-chance-randomness}
Antony Eagle.
\newblock {Chance versus Randomness}.
\newblock In \emph{The {Stanford} Encyclopedia of Philosophy}. Metaphysics Research Lab, Stanford University, {S}pring 2021 edition, 2021.

\bibitem[Elga(2010)]{elga2010subjective}
Adam Elga.
\newblock Subjective probabilities should be sharp.
\newblock \emph{Philosophers' Imprint}, 10\penalty0 (5):\penalty0 1--11, 2010.

\bibitem[Ellsberg(1961)]{ellsberg1961risk}
Daniel Ellsberg.
\newblock Risk, ambiguity, and the {S}avage axioms.
\newblock \emph{The Quarterly Journal of Economics}, 75\penalty0 (4):\penalty0 643--669, 1961.

\bibitem[Epstein \& Schneider(2003)Epstein and Schneider]{epstein2003iid}
Larry~G. Epstein and Martin Schneider.
\newblock {IID}: independently and indistinguishably distributed.
\newblock \emph{Journal of Economic Theory}, 113\penalty0 (1):\penalty0 32--50, 2003.

\bibitem[Feller(1991)]{feller1991introduction}
William Feller.
\newblock \emph{An introduction to probability theory and its applications, Volume 2, 2nd edition}.
\newblock John Wiley \& Sons, 1991.

\bibitem[Fierens(2003)]{fierens2003towardsphdthesis}
Pablo~I. Fierens.
\newblock \emph{Towards a chaotic probability model for frequentist probability}.
\newblock PhD thesis, Cornell University, 2003.

\bibitem[Fierens(2009)]{fierens2009extension}
Pablo~I. Fierens.
\newblock An extension of chaotic probability models to real-valued variables.
\newblock \emph{International Journal of Approximate Reasoning}, 50\penalty0 (4):\penalty0 627--641, 2009.

\bibitem[Fierens \& Fine(2001)Fierens and Fine]{fierens2001towards}
Pablo~I. Fierens and Terrence~L Fine.
\newblock Towards a frequentist interpretation of sets of measures.
\newblock In \emph{International Symposium on Imprecise Probabilities: Theories and Applications}, pp.\  179--187, 2001.

\bibitem[Fierens \& Fine(2003)Fierens and Fine]{fierens2003towards}
Pablo~I. Fierens and Terrence~L. Fine.
\newblock Towards a chaotic probability model for frequentist probability: The univariate case.
\newblock In \emph{International Symposium on Imprecise Probabilities: Theories and Applications}, pp.\  245--259, 2003.

\bibitem[Fierens et~al.(2009)Fierens, R{\^e}go, and Fine]{fierens2009frequentist}
Pablo~I. Fierens, Leonardo~C. R{\^e}go, and Terrence~L. Fine.
\newblock A frequentist understanding of sets of measures.
\newblock \emph{Journal of Statistical Planning and Inference}, 139\penalty0 (6):\penalty0 1879--1892, 2009.

\bibitem[Fine(1970)]{fine1970apparent}
Terrence~L. Fine.
\newblock On the apparent convergence of relative frequency and its implications.
\newblock \emph{IEEE Transactions on Information Theory}, 16\penalty0 (3):\penalty0 251--257, 1970.

\bibitem[Fine(1988)]{fine1988lower}
Terrence~L. Fine.
\newblock Lower probability models for uncertainty and nondeterministic processes.
\newblock \emph{Journal of Statistical Planning and Inference}, 20\penalty0 (3):\penalty0 389--411, 1988.

\bibitem[Franzke et~al.(2020)Franzke, Barbosa, Blender, Fredriksen, Laepple, Lambert, Nilsen, Rypdal, Rypdal, Scotto, et~al.]{franzke2020structure}
Christian~L.E. Franzke, Susana Barbosa, Richard Blender, Hege-Beate Fredriksen, Thomas Laepple, Fabrice Lambert, Tine Nilsen, Kristoffer Rypdal, Martin Rypdal, Manuel~G. Scotto, et~al.
\newblock The structure of climate variability across scales.
\newblock \emph{Reviews of Geophysics}, 58\penalty0 (2):\penalty0 1--44, 2020.
\newblock Paper number e2019RG000657.

\bibitem[Fr{\"o}hlich \& Williamson(2024)Fr{\"o}hlich and Williamson]{frohlich2022risk}
Christian Fr{\"o}hlich and Robert~C. Williamson.
\newblock Risk measures and upper probabilities: Coherence and stratification.
\newblock \emph{Journal of Machine Learning Research (to appear)}, 2024.

\bibitem[Fr{\"o}hlich et~al.(2024)Fr{\"o}hlich, Derr, and Williamson]{frohlich2024strictly}
Christian Fr{\"o}hlich, Rabanus Derr, and Robert~C. Williamson.
\newblock Strictly frequentist imprecise probability.
\newblock \emph{International Journal of Approximate Reasoning}, 168, 2024.
\newblock Paper number 109148.

\bibitem[Galvan(2006)]{galvan2006bohmian}
Bruno Galvan.
\newblock Bohmian mechanics and typicality without probability.
\newblock \emph{arXiv preprint quant-ph/0605162}, 2006.

\bibitem[Gilboa \& Schmeidler(1989)Gilboa and Schmeidler]{gilboa1989maxmin}
Itzhak Gilboa and David Schmeidler.
\newblock Maxmin expected utility with non-unique prior.
\newblock \emph{Journal of {M}athematical {E}conomics}, 18\penalty0 (2):\penalty0 141--153, 1989.

\bibitem[Gneiting \& Raftery(2007)Gneiting and Raftery]{gneiting2007strictly}
Tilmann Gneiting and Adrian~E. Raftery.
\newblock Strictly proper scoring rules, prediction, and estimation.
\newblock \emph{Journal of the American Statistical Association}, 102\penalty0 (477):\penalty0 359--378, 2007.

\bibitem[Gorban(2017)]{gorban2017statistical}
Igor~I. Gorban.
\newblock \emph{The Statistical Stability Phenomenon}.
\newblock Springer, 2017.

\bibitem[Grize \& Fine(1987)Grize and Fine]{grize1987continuous}
Yves~L. Grize and Terrence~L. Fine.
\newblock Continuous lower probability-based models for stationary processes with bounded and divergent time averages.
\newblock \emph{The Annals of Probability}, 15\penalty0 (2):\penalty0 783--803, 1987.

\bibitem[Haghtalab et~al.(2022)Haghtalab, Jordan, and Zhao]{haghtalab2022demand}
Nika Haghtalab, Michael Jordan, and Eric Zhao.
\newblock On-demand sampling: Learning optimally from multiple distributions.
\newblock In \emph{Advances in Neural Information Processing Systems}, volume~35, pp.\  406--419. 2022.

\bibitem[Hu et~al.(2016)Hu, Chen, and Wu]{hu2016general}
Feng Hu, Zengjing Chen, and Panyu Wu.
\newblock A general strong law of large numbers for non-additive probabilities and its applications.
\newblock \emph{Statistics}, 50\penalty0 (4):\penalty0 733--749, 2016.

\bibitem[Huber(1981)]{huber1981robust}
Peter~J. Huber.
\newblock \emph{Robust Statistics}.
\newblock John Wiley \& Sons, Inc., 1981.

\bibitem[H{\"u}llermeier et~al.(2022)H{\"u}llermeier, Destercke, and Shaker]{hullermeier2022quantification}
Eyke H{\"u}llermeier, S{\'e}bastien Destercke, and Mohammad~Hossein Shaker.
\newblock Quantification of credal uncertainty in machine learning: A critical analysis and empirical comparison.
\newblock In \emph{Uncertainty in Artificial Intelligence}, pp.\  548--557. PMLR, 2022.

\bibitem[Iacovissi et~al.(2023)Iacovissi, Lu, and Williamson]{iacovissi2023general}
Laura Iacovissi, Nan Lu, and Robert~C Williamson.
\newblock A general framework for learning under corruption: Label noise, attribute noise, and beyond.
\newblock \emph{arXiv preprint arXiv:2307.08643}, 2023.

\bibitem[Ibragimov et~al.(2015)Ibragimov, Ibragimov, and Walden]{ibragimov2015heavy}
Marat Ibragimov, Rustam Ibragimov, and Johan Walden.
\newblock \emph{Heavy-tailed distributions and robustness in economics and finance}.
\newblock Springer, 2015.

\bibitem[Ivanenko(2010)]{ivanenkobook}
Victor~I. Ivanenko.
\newblock \emph{Decision Systems and Nonstochastic Randomness}.
\newblock Springer, 2010.

\bibitem[Ivanenko \& Pasichnichenko(2017)Ivanenko and Pasichnichenko]{ivanenko2017expected}
Victor~I. Ivanenko and Illia Pasichnichenko.
\newblock Expected utility for nonstochastic risk.
\newblock \emph{Mathematical Social Sciences}, 86:\penalty0 18--22, 2017.

\bibitem[Joshi(1983)]{joshi1983introduction}
Kapil~D. Joshi.
\newblock \emph{Introduction to general topology}.
\newblock New Age International, 1983.

\bibitem[Klenke(2013)]{klenke2013probability}
Achim Klenke.
\newblock \emph{Probability theory: a comprehensive course}.
\newblock Springer, 2013.

\bibitem[Konek(2019)]{konek2019ip}
Jason Konek.
\newblock Ip scoring rules: foundations and applications.
\newblock In \emph{International Symposium on Imprecise Probabilities: Theories and Applications}, pp.\  256--264. PMLR, 2019.

\bibitem[Kuhn et~al.(2019)Kuhn, Esfahani, Nguyen, and Shafieezadeh-Abadeh]{kuhn2019wasserstein}
Daniel Kuhn, Peyman~Mohajerin Esfahani, Viet~Anh Nguyen, and Soroosh Shafieezadeh-Abadeh.
\newblock Wasserstein distributionally robust optimization: Theory and applications in machine learning.
\newblock In \emph{Operations Research \& Management Science in the Age of Analytics}, pp.\  130--166. Informs, 2019.

\bibitem[Kumar \& Fine(1985)Kumar and Fine]{kumar1985stationary}
Anurag Kumar and Terrence~L. Fine.
\newblock Stationary lower probabilities and unstable averages.
\newblock \emph{Zeitschrift f{\"u}r Wahrscheinlichkeitstheorie und verwandte Gebiete}, 69\penalty0 (1):\penalty0 1--17, 1985.

\bibitem[La~Caze(2016)]{lacazefrequentism}
Adam La~Caze.
\newblock {Frequentism}.
\newblock In Alan Hájek and Christopher Hitchcock (eds.), \emph{{The Oxford Handbook of Probability and Philosophy}}. Oxford University Press, 2016.

\bibitem[Lovejoy(2015)]{lovejoy2015voyage}
S.~Lovejoy.
\newblock A voyage through scales, a missing quadrillion and why the climate is not what you expect.
\newblock \emph{Climate Dynamics}, 44:\penalty0 3187--3210, 2015.

\bibitem[Maccheroni \& Marinacci(2005)Maccheroni and Marinacci]{maccheroni2005strong}
Fabio Maccheroni and Massimo Marinacci.
\newblock A strong law of large numbers for capacities.
\newblock \emph{The Annals of Probability}, 33\penalty0 (3):\penalty0 1171--1178, 2005.

\bibitem[Marinacci(1999)]{marinacci1999limit}
Massimo Marinacci.
\newblock Limit laws for non-additive probabilities and their frequentist interpretation.
\newblock \emph{Journal of Economic Theory}, 84\penalty0 (2):\penalty0 145--195, 1999.

\bibitem[Mayo-Wilson \& Wheeler(2016)Mayo-Wilson and Wheeler]{mayo2016scoring}
Conor Mayo-Wilson and Gregory Wheeler.
\newblock Scoring imprecise credences.
\newblock \emph{Philosophy and Phenomenological Research}, 93\penalty0 (1):\penalty0 55--78, 2016.

\bibitem[McMahan et~al.(2017)McMahan, Moore, Ramage, Hampson, and y~Arcas]{mcmahan2017communication}
Brendan McMahan, Eider Moore, Daniel Ramage, Seth Hampson, and Blaise~Aguera y~Arcas.
\newblock Communication-efficient learning of deep networks from decentralized data.
\newblock In \emph{International Conference on Artificial Intelligence and Statistics}, volume~54, pp.\  1273--1282. PMLR, 2017.

\bibitem[Miranda \& Cooman(2014)Miranda and Cooman]{introtoiplowerprev}
Enrique Miranda and Gert~de Cooman.
\newblock Lower previsions.
\newblock In \emph{Introduction to Imprecise Probabilities}, pp.\  28--55. John Wiley \& Sons, Ltd, 2014.

\bibitem[Nair et~al.(2022)Nair, Wierman, and Zwart]{nair2022fundamentals}
Jayakrishnan Nair, Adam Wierman, and Bert Zwart.
\newblock \emph{The fundamentals of heavy tails: Properties, emergence, and estimation}.
\newblock Cambridge University Press, 2022.

\bibitem[Papamarcou \& Fine(1991{\natexlab{a}})Papamarcou and Fine]{papamarcou1991stationarity}
Adrian Papamarcou and Terrence~L. Fine.
\newblock Stationarity and almost sure divergence of time averages in interval-valued probability.
\newblock \emph{Journal of Theoretical Probability}, 4\penalty0 (2):\penalty0 239--260, 1991{\natexlab{a}}.

\bibitem[Papamarcou \& Fine(1991{\natexlab{b}})Papamarcou and Fine]{papamarcou1991unstable}
Adrian Papamarcou and Terrence~L. Fine.
\newblock Unstable collectives and envelopes of probability measures.
\newblock \emph{The Annals of Probability}, 19\penalty0 (2):\penalty0 893--906, 1991{\natexlab{b}}.

\bibitem[Peng \& Qi(2018)Peng and Qi]{peng2017inference}
Liang Peng and Yongcheng Qi.
\newblock \emph{Inference for heavy-tailed data: applications in insurance and finance}.
\newblock Academic press, 2018.

\bibitem[Peng(2007)]{peng2007law}
Shige Peng.
\newblock Law of large numbers and central limit theorem under nonlinear expectations.
\newblock \emph{arXiv preprint math/0702358}, 2007.

\bibitem[Peng(2019)]{peng2019nonlinear}
Shige Peng.
\newblock \emph{Nonlinear expectations and stochastic calculus under uncertainty: with robust CLT and G-Brownian motion}.
\newblock Springer Nature, 2019.

\bibitem[Perdomo et~al.(2020)Perdomo, Zrnic, Mendler-D{\"u}nner, and Hardt]{perdomo2020performative}
Juan Perdomo, Tijana Zrnic, Celestine Mendler-D{\"u}nner, and Moritz Hardt.
\newblock Performative prediction.
\newblock In \emph{International Conference on Machine Learning}, volume 119, pp.\  7599--7609. PMLR, 2020.

\bibitem[Persiau et~al.(2021)Persiau, De~Bock, and De~Cooman]{persiau2021remarkable}
Floris Persiau, Jasper De~Bock, and Gert De~Cooman.
\newblock A remarkable equivalence between non-stationary precise and stationary imprecise uncertainty models in computable randomness.
\newblock In \emph{International Symposium on Imprecise Probability: Theories and Applications}, volume 147, pp.\  244--253. PMLR, 2021.

\bibitem[Persiau et~al.(2022)Persiau, {De Bock}, and {de Cooman}]{persiau2022on}
Floris Persiau, Jasper {De Bock}, and Gert {de Cooman}.
\newblock On the (dis)similarities between stationary imprecise and non-stationary precise uncertainty models in algorithmic randomness.
\newblock \emph{International Journal of Approximate Reasoning}, 151:\penalty0 272--291, 2022.

\bibitem[Punzo et~al.(2018)Punzo, Bagnato, and Maruotti]{PUNZO201895}
Antonio Punzo, Luca Bagnato, and Antonello Maruotti.
\newblock Compound unimodal distributions for insurance losses.
\newblock \emph{Insurance: Mathematics and Economics}, 81:\penalty0 95--107, 2018.

\bibitem[R{\'e}bill{\'e}(2009)]{rebille2009law}
Yann R{\'e}bill{\'e}.
\newblock Law of large numbers for non-additive measures.
\newblock \emph{Journal of Mathematical Analysis and Applications}, 352\penalty0 (2):\penalty0 872--879, 2009.

\bibitem[R{\^e}go \& Fine(2005)R{\^e}go and Fine]{rego2005estimation}
Leandro~C. R{\^e}go and Terrence~L. Fine.
\newblock Estimation of chaotic probabilities.
\newblock In \emph{International Symposium on Imprecise Probabilities: Theories and Applications}, volume~5, pp.\  297--305, 2005.

\bibitem[Resnick(2007)]{resnick2007heavy}
Sidney~I. Resnick.
\newblock \emph{Heavy-tail phenomena: probabilistic and statistical modeling}.
\newblock Springer Science \& Business Media, 2007.

\bibitem[Rockafellar \& Royset(2015)Rockafellar and Royset]{rockafellar2015measures}
R.~Tyrrell Rockafellar and Johannes~O Royset.
\newblock Measures of residual risk with connections to regression, risk tracking, surrogate models, and ambiguity.
\newblock \emph{SIAM Journal on Optimization}, 25\penalty0 (2):\penalty0 1179--1208, 2015.

\bibitem[Sadrolhefazi \& Fine(1994)Sadrolhefazi and Fine]{sadrolhefazi1994finite}
Amir Sadrolhefazi and Terrence~L. Fine.
\newblock Finite-dimensional distributions and tail behavior in stationary interval-valued probability models.
\newblock \emph{The Annals of Statistics}, 22\penalty0 (4):\penalty0 1840--1870, 1994.

\bibitem[Schoenfield(2017)]{schoenfield2017accuracy}
Miriam Schoenfield.
\newblock The accuracy and rationality of imprecise credences.
\newblock \emph{No{\^u}s}, 51\penalty0 (4):\penalty0 667--685, 2017.

\bibitem[Seidenfeld et~al.(2012)Seidenfeld, Schervish, and Kadane]{seidenfeld2012forecasting}
Teddy Seidenfeld, Mark~J Schervish, and Joseph~B Kadane.
\newblock Forecasting with imprecise probabilities.
\newblock \emph{International Journal of Approximate Reasoning}, 53\penalty0 (8):\penalty0 1248--1261, 2012.

\bibitem[Shafer \& Vovk(2001)Shafer and Vovk]{shafervovk2001probabilityfinance}
Glenn Shafer and Vladimir Vovk.
\newblock \emph{Probability and finance: it's only a game!}
\newblock John Wiley \& Sons, 2001.

\bibitem[Spreij()]{spreij2012measure}
Peter~J.C. Spreij.
\newblock Measure theoretic probability.
\newblock \emph{Lecture Notes}.

\bibitem[Taleb(2007)]{taleb2007black}
Nassim~Nicholas Taleb.
\newblock \emph{The black swan: The impact of the highly improbable}.
\newblock Random House, 2007.

\bibitem[Ter{\'a}n(2014)]{teran2014laws}
Pedro Ter{\'a}n.
\newblock Laws of large numbers without additivity.
\newblock \emph{Transactions of the American Mathematical Society}, 366\penalty0 (10):\penalty0 5431--5451, 2014.

\bibitem[Trautmann \& Van De~Kuilen(2015)Trautmann and Van De~Kuilen]{trautmann2015ambiguity}
Stefan~T. Trautmann and Gijs Van De~Kuilen.
\newblock Ambiguity attitudes.
\newblock In \emph{The Wiley Blackwell handbook of judgment and decision making}, pp.\  89--116. John Wiley \& Sons, Ltd., 2015.

\bibitem[T’Joens et~al.(2019)T’Joens, Krak, Bock, and Cooman]{t2019recursive}
Natan T’Joens, Thomas Krak, Jasper~De Bock, and Gert~de Cooman.
\newblock A recursive algorithm for computing inferences in imprecise markov chains.
\newblock In \emph{Symbolic and Quantitative Approaches to Reasoning with Uncertainty: 15th European Conference, ECSQARU 2019, Belgrade, Serbia, September 18-20, 2019, Proceedings 15}, pp.\  455--465. Springer, 2019.

\bibitem[Vershynin(2020)]{vershynin2018high}
Roman Vershynin.
\newblock \emph{High-dimensional probability: An introduction with applications in data science}.
\newblock 2020.
\newblock URL \url{https://www.math.uci.edu/~rvershyn/papers/HDP-book/HDP-book.html}.
\newblock Version of 2020-06-09.

\bibitem[von Mises(1919)]{mises1919grundlagen}
Richard von Mises.
\newblock {Grundlagen der Wahrscheinlichkeitsrechnung}.
\newblock \emph{Mathematische Zeitschrift}, 5\penalty0 (1):\penalty0 52--99, 1919.

\bibitem[Vovk et~al.(2005)Vovk, Takemura, and Shafer]{vovk2005defensive}
Vladimir Vovk, Akimichi Takemura, and Glenn Shafer.
\newblock Defensive forecasting.
\newblock In \emph{International Workshop on Artificial Intelligence and Statistics}, pp.\  365--372. PMLR, 2005.

\bibitem[Walley(1991)]{walley1991statistical}
Peter Walley.
\newblock \emph{Statistical reasoning with imprecise probabilities}.
\newblock Chapman-Hall, 1991.

\bibitem[Walley \& Fine(1982)Walley and Fine]{walley1982towards}
Peter Walley and Terrence~L. Fine.
\newblock Towards a frequentist theory of upper and lower probability.
\newblock \emph{The Annals of Statistics}, 10\penalty0 (3):\penalty0 741--761, 1982.

\bibitem[Williams(1991)]{williams1991probability}
David Williams.
\newblock \emph{Probability with martingales}.
\newblock Cambridge {U}niversity {P}ress, 1991.

\bibitem[Williamson \& Menon(2019)Williamson and Menon]{williamson2019fairness}
Robert~C. Williamson and Aditya Menon.
\newblock Fairness risk measures.
\newblock In \emph{International Conference on Machine Learning}, volume~97, pp.\  6786--6797. PMLR, 2019.

\bibitem[Zhang et~al.(2021)Zhang, Xie, Bai, Yu, Li, and Gao]{zhang2021survey}
Chen Zhang, Yu~Xie, Hang Bai, Bin Yu, Weihong Li, and Yuan Gao.
\newblock A survey on federated learning.
\newblock \emph{Knowledge-Based Systems}, 216, 2021.
\newblock Paper number 106775.

\bibitem[Zhang et~al.(2024)Zhang, Tang, and Xiong]{zhang2024conditional}
Jiaqi Zhang, Yanyan Tang, and Jie Xiong.
\newblock Conditional strong law of large numbers under g-expectations.
\newblock \emph{Symmetry}, 16\penalty0 (3):\penalty0 272, 2024.

\bibitem[Zhang \& Lan(2020)Zhang and Lan]{zhang2020strong}
Ning Zhang and Yuting Lan.
\newblock A strong law of large numbers for independent random variables under non-additive probabilities.
\newblock \emph{Communications in Statistics --- Theory and Methods}, 49\penalty0 (21):\penalty0 5252--5272, 2020.

\bibitem[Zhao \& Ermon(2021)Zhao and Ermon]{zhao2021right}
Shengjia Zhao and Stefano Ermon.
\newblock Right decisions from wrong predictions: A mechanism design alternative to individual calibration.
\newblock In \emph{International Conference on Artificial Intelligence and Statistics}, pp.\  2683--2691. PMLR, 2021.

\end{thebibliography}
\bibliographystyle{tmlr}

\clearpage

\appendix
\section{Appendix}

\subsection{Proof of Proposition~\ref{prop:galvantypicality}}
\label{app:galvantypicalityproof}
\begin{proof}
A typicality space in the sense of \citet{galvan2006bohmian} is a triplet $(\Lambda,\mathcal{A},d)$, where $\Lambda$ is a ``sample space'' (not further defined), $\mathcal{A} \subseteq 2^\Lambda$ is a set system closed under complementation and which contains the emptyset, and $d$ is a \textit{typicality distance}. 
We show that $d(A,B) \coloneqq \Pup(A \triangle B)$, $A,B \in \mathcal{A}$, satisfies the axioms of a typicality distance as in \citep{galvan2006bohmian}, which are as follows:
\begin{enumerate}[nolistsep,start=1,label=\textbf{T\arabic*.}, ref=T\arabic*]
\item \label{item:T1} $A \subseteq A' \Rightarrow d(A,\emptyset)\leq d(A',\emptyset)$. 
\item \label{item:T2} $d(\Lambda,\emptyset)=1$.
\item \label{item:T3} $d(A,\Lambda)=d(A^c,\emptyset)$.
\item \label{item:T4} If $A_1 \cap A_2 \in \mathcal{A}$, then $d(A_1\cap A_2,B) \leq d(A_1,B) + d(A_2,B)$.
\item \label{item:T5} If $A_1\cap A_2$, $B_1 \cap B_2 \in \mathcal{A}$, then $d(A_1 \cap A_2, B_1 \cap B_2) \leq d(A_1,B_1) + d(A_2,B_2)$.
\end{enumerate}
We note that if $P$ is a probability measure, then $\tilde{d}(A,B)\coloneqq P(A \triangle B)$ satisfies these axioms. Since $\Pup$ is sub-additive in the sense that $\Pup(A\cup B)\leq \Pup(A) +\Pup(B)$, we find that the relevant properties are preserved.

\ref{item:T1} follows from monotonicity of $\Pup$; note that $d(A,\emptyset)=\Pup(A)$. \ref{item:T2} follows since $\Pup(\Lambda)=1$. For~\ref{item:T3}, note that $d(A,\Lambda)=\Pup(A \triangle \Lambda)=\Pup(A^c)$ and $d(A^c,\emptyset)=\Pup(A^c)$. To see that \ref{item:T4} and \ref{item:T5} hold, it is easiest so draw Venn diagrams. For completeness we provide the computations.

For \ref{item:T4}, observe that $d(A_1 \cap A_2,B)=\Pup(\RomanNumeralCaps{1})$ with
\begin{align}
    \RomanNumeralCaps{1} \coloneqq &(A_1 \cap A_2) \triangle B = (A_1 \cap A_2 \cap B^C) \cup (A_1^C \cup A_2^C) \cap B\\
    &\subseteq (A_1^C \cap B) \cup (A_1 \cap B^C) \cup (A_2^C \cap B) \cup (A_2 \cap B^C) \eqqcolon \RomanNumeralCaps{2}.
\end{align}
From monotonicity of $\Pup$ we get that $\Pup(\RomanNumeralCaps{1}) \leq \Pup(\RomanNumeralCaps{2})$, and from subadditivity 
\[
\Pup(\RomanNumeralCaps{2}) \leq \Pup((A_1^C \cap B) \cup (A_1 \cap B^C)) + \Pup((A_2^C \cap B) \cup (A_2 \cap B^C)) = d(A_1,B) + d(A_2,B).
\]
For~\ref{item:T5}, $d(A_1 \cap A_2, B_1 \cap B_2)=\Pup(\RomanNumeralCaps{1})$, where $\RomanNumeralCaps{1} \coloneqq \left(A_1 \cap A_2 \cap (B_1\cap B_2)^c\right) \cup \left((A_1\cap A_2)^c \cap (B_1 \cap B_2) \right)$. Now,
\begin{align}
    \RomanNumeralCaps{1} \subseteq (A_1 \cap B_1^c) \cup (A_1^c \cap B_1) \cup (A_2 \cap B_2^c) \cup (A_2^c \cap B_2) \eqqcolon \RomanNumeralCaps{2},
\end{align}
because \[
(A_1^c \cup A_2^c) \cap B_1 \cap B_2 = (A_1^c \cap B_1 \cap B_2) \cup (A_2^c \cap B_1 \cap B_2) \subseteq (A_1^c \cap B_1) \cup (A_2^c \cap B_2)
\]
and 
\[
A_1 \cap A_2 \cap (B_1^c \cup B_2^c) = (A_1 \cap A_2 \cap B_1^c) \cup (A_1 \cap A_2 \cap B_2^c) \subseteq (A_1 \cap B_1^c) \cup (A_2 \cap B_2^c).
\]
From $\RomanNumeralCaps{1} \subseteq \RomanNumeralCaps{2}$, we conclude by monotonicity and subadditivity that 
\[\Pup(\RomanNumeralCaps{1}) \leq \Pup\left((A_1 \cap B_1^c) \cup (A_1^c \cap B_1)\right) + \Pup\left((A_2 \cap B_2^c) \cup (A_2^c \cap B_2)\right) = d(A_1,B_1) + d(A_2,B_2).
\]

Based on a typicality distance, \citet{galvan2006bohmian} defines a \textit{relative typicality measure} as 
\[
T_r(A|B) \coloneqq \inf_{A' \subseteq A, B' \supseteq B} \frac{d(A',B')}{d(B,\emptyset)}, \quad A,B \in \mathcal{A},
\]
and the \textit{absolute typicality measure} as 
\[
T_a(A) \coloneqq T_r(A|\Lambda) = d(A^c,\emptyset), \quad A \in \mathcal{A}.
\]
In our case, $T_a(A) = \Pup(A^c \triangle \emptyset)=\Pup(A^c) = 1-\Plow(A)$. 
\end{proof}

\subsection{Proof of the Main Theorem~\ref{theorem:maintheoremffivanenko}}
\label{sec:proofofmaintheorem}
Throughout this section, we work on the measure space $([0,1), \mathcal{B}([0,1)), \lambda)$.
We here prove our main theorem, restated for convenience:
\begin{theorem}

Let $\mathcal{Z}=[0,1)$ and $\emptyset \neq \mathcal{M} \subseteq \Delta^k$. Then there exists a sequence of probability measures $m_1,m_2,.. \in \mM$ with independent product
    \begin{equation}
        \label{eq:independentproduct}
        \lambda\left\{z \in \mathcal{Z} : (\omega_1,\omega_2,..,\omega_n) \subset G(z)\}\right) \coloneqq \prod_{i=1}^n m_i(\{\omega_i\}), \quad \forall n \in \mathbb{N},
    \end{equation}
 so that it holds
    \begin{equation}
    \label{eq:divevent}
    \lambda\left\{z \in \mathcal{Z}:  \CP\left(n \mapsto r_n^{G(z)}\right) =\CP\left(n \mapsto \frac{1}{n} \sumin m_i \right) = \mathcal{N}\right\} = 1
    \end{equation}
if and only if $\emptyset \neq \mathcal{N} \subseteq \cobar(\mathcal{M})$ is a closed connected subset of the closed convex hull of $\mM$. 

We can in addition demand that any $m \in \mM$ is a cluster point of the sequence of measures $m_i$ (with respect to the Euclidean metric).
\end{theorem}

We build this up from some helper propositions.
\begin{proposition}
\label{prop:samecpastheoretical}
    Let $p^\infty=(p_1,p_2,..)$ be a sequence of probability measures in $\Delta^k$, with independent product
    \[
\lambda\left\{z \in \mathcal{Z} : (\omega_1,\omega_2,..,\omega_n) \subset G(z)\}\right) \coloneqq \prod_{i=1}^n p_i(\{\omega_i\}), \quad \forall n \in \mathbb{N}.
\] 
Then
 \[
    \lambda\left\{z \in \mathcal{Z}:  \CP\left(n \mapsto r_n^{G(z)}\right) = \CP\left(n \mapsto \frac{1}{n} \sum_{i=1}^n p_i\right)\right\} = 1.
    \]
\end{proposition}
In words, almost surely the cluster points of relative frequencies will coincide exactly with the cluster points of the sequence of measures (with respect to the Euclidean topology on $\Delta^k$).

This follows by applying a variant of the strong law of large numbers, specialized to our context:
\begin{theorem}[{\citep[Section 10.7]{feller1991introduction}}]
\label{theorem:kolmogorovsstronglaw}
Let $X_1,X_2,..$ be independent but not necessarily identically distributed random variables. Define $\bar{X}_n \coloneqq \frac{1}{n} \sum_{i=1}^n X_i$. If they satisfy
\[
Var(X_k)<\infty \; \forall k \in \mathbb{N} \;\; \text{ and } \;\; \sum_{k=1}^\infty \frac{Var(X_k)}{k^2} < \infty,
\]
then it holds
\[
\lambda\left\{\lim_{n \rightarrow \infty} \bar{X}_n - \mathbb{E}[\bar{X}_n] = 0\right\} = 1.
\]
In words, the sample average converges almost surely to its expectation.    
\end{theorem}
\begin{remark}\normalfont
    The statement can be extended to vector-valued random variables, i.e $X_i \in \mR^d$, when we check the conditions for each component. Then we get $d$-many almost sure events (almost sure convergence in each dimension), and intersecting finitely many almost sure events yields an almost sure event.
\end{remark}
We now prove Proposition~\ref{prop:samecpastheoretical}.
\begin{proof}
Under our independent product, the projections $W_i : [0,1) \to \Omega$ are measurable.
    Consider the sequence of $\{0,1\}^k$-valued random variables $\chi_1,\chi_2,..$ where 
    \[
    \chi_i(z) \coloneqq \left(\chi_{\{\omega^1\}}(W_i(z)), .., \chi_{\{\omega^k\}}(W_i(z))\right)^\intercal.
    \]
    This sequence satisfies the condition of Theorem~\ref{theorem:kolmogorovsstronglaw}:

     First, we check that each component $\chi_i^j(z)$, $j=1..k$, of $\chi_i(z)$ has finite variance; but its variance is simply $\operatorname{Var}(\chi_i^j) = p_i(\omega^j)(1- p_i(\omega^j))\leq 0.25 < \infty$ (Bernoulli variance). But then it is easy to see that for any component $j$
    \[
    \sum_{i=1}^\infty \frac{\operatorname{Var}(\chi_i^j)}{i^2}  < \infty.
    \]
    
    Applying Theorem~\ref{theorem:kolmogorovsstronglaw} therefore gives
    \begin{align}
    &\lambda\left\{\lim_{n\to \infty} \frac{1}{n} \sum_{i=1}^n \chi_i(z) - \mathbb{E}\left[\frac{1}{n} \sum_{i=1}^n \chi_i(z)\right] = 0\right\} = 1\\
    \Leftrightarrow \; &\lambda\left\{\lim_{n\to \infty} r_n^{G(z)} - \mathbb{E}\left[r_n^{G(z)}\right] = 0\right\} = 1\\
    \Leftrightarrow \; &\lambda\left\{\lim_{n\to \infty} r_n^{G(z)} - \frac{1}{n} \sum_{i=1}^n p_i = 0\right\} = 1.
    \end{align}
    It remains to observe that the last condition implies
    \[
    \lambda\left\{z \in \mathcal{Z}:  \CP\left(n \mapsto r_n^{G(z)}\right) = \CP\left(n \mapsto \frac{1}{n} \sum_{i=1}^n p_i\right)\right\} = 1.
    \]
    But it is clear that if for two $\mR^k$ valued sequences $a^\infty=(a_1,..),b^\infty=(b_1,..)$, we have $\limn a_n - b_n = 0$, then $\CP(a^\infty)=\CP(b^\infty)$. We briefly demonstrate this. Assume $c$ is a cluster point of $a^\infty$, which is equivalent to 
    \[
\forall \varepsilon>0 : \forall n_0 \in \mathbb{N} : \exists n \geq n_0 : a_{n} \in B_\varepsilon(c).
\]
We want to show that then $c$ is also a cluster point of $b^\infty$, equivalent to $\forall \varepsilon'>0 : \forall n_0' \in \mathbb{N} : \exists n \geq n_0': b_{n} \in B_{\varepsilon'}(c)$. But since  $\limn a_n - b_n = 0$, we can find $\varepsilon+\kappa<\varepsilon'$ so that $\exists n_0'' \in \mathbb{N} : \forall n \geq n_0'' : d(a_n,b_n)<\kappa$. Setting $n_0 \coloneqq n_0' \coloneqq n_0''$ then gives us some $n \geq n_0$ for which $a_n \in B_\varepsilon(c)$ and $d(a_n,b_n)<\kappa$, implying that $d(c,b_n)<\varepsilon+\kappa<\varepsilon'$.
\end{proof}

Thus, in light of Proposition~\ref{prop:samecpastheoretical}, in order to prove Theorem~\ref{theorem:maintheoremffivanenko}, it remains to show the following.
\begin{proposition}
\label{prop:mainproptheoretical}
    Let $\emptyset \neq \mathcal{M} \subseteq \Delta^k$. Then there exists a sequence of measures $m^\infty=(m_1,m_2,..)$, $m_i \in \mM$, so that
    \[
    \CP\left(n \mapsto \frac{1}{n} \sum_{i=1}^n m_i\right) = \mathcal{N},
    \]
    if and only if $\emptyset \neq \mathcal{N} \subseteq \cobar(\mathcal{M})$ is a closed connected subset of the closed convex hull of $\mM$. 
    We can in addition demand that any point $m \in \mM$ is a cluster point of the sequence of measures $m_i$ (with respect to the Euclidean metric).
\end{proposition}

For the forward direction, the key proposition of interest is the following.
\begin{proposition}
\label{prop:forwarddir}
    Let $\emptyset \neq \mathcal{M} \subseteq \Delta^k$ and $\emptyset \neq \mathcal{N} \subseteq \cobar(\mathcal{M})$ a closed connected subset of the closed convex hull of $\mM$. Then there exists a sequence of measures $m^\infty=(m_1,m_2,..)$, $m_i \in \mM$, so that
    \[
    \CP\left(n \mapsto \frac{1}{n} \sum_{i=1}^n m_i\right) = \mathcal{N}.
    \] 
    In addition, it can be guaranteed that any point $n \in \mathcal{N}$ is a cluster point of the sequence of measures $m_i$.
\end{proposition}
We defer its proof to later.
The forward direction is concluded by showing in addition:
\begin{proposition}
    Let $\emptyset \neq \mathcal{M} \subseteq \Delta^k$ and $m^\infty=(m_1,m_2,..)$ a sequence of measures from $\mM$. Then there exists a sequence of measure $\tilde{m}^\infty=(\tilde{m}_1,\tilde{m}_2,..)$ from $\mM$ so that any $m \in \mM$ is a cluster point of $\tilde{m}^\infty$ and 
    \[
    \CP\left(n \mapsto \frac{1}{n} \sum_{i=1}^n m_i\right) = \CP\left(n \mapsto \frac{1}{n} \sum_{i=1}^n \tilde{m}_i\right).
    \]
\end{proposition}
\begin{proof}
 Note that $\cobar(\mM)$ itself is closed and connected, hence satisfies the assumptions of Proposition~\ref{prop:forwarddir}. Thus we can obtain a sequence $m^{\infty,2}$ for which every measure in $\mM$ is a cluster point (indeed any measure from $\cobar(\mM)$). Now construct a selection rule $S: \mathbb{N} \to \{0,1\}$, with ``rate going to zero'', meaning $\limn \frac{\sum_{i=1}^n S(n)}{n} = 0$. We use this selection rule to ``interleave'' the sequences $m^{\infty,2}$ and $m^\infty$ in the following way. The new sequence of measures is
 \[
 \tilde{m}_i \coloneqq \begin{cases}
     m_{j}^{2} : j \coloneqq \sum_{j=1}^i S(j) & \text{ if } S(i)=1\\
     m_{j} : j \coloneqq i-\sum_{j=1}^i S(j) & \text{ if } S(i)=0.\\
 \end{cases}
 \]
From this we obtain
  \[
    \frac{1}{n}\sum_{i=1}^n  \tilde{m}_i =  \frac{1}{n} \sum_{i=1}^{\sum_{j=1}^n S(j)} m_{i}^{2}  +  \frac{1}{n} \sum_{i=1}^{n-\sum_{j=1}^n S(j)} m_j .
    \]
   But since $\limn \frac{\sum_{i=1}^n S(n)}{n} = 0$ and all measures in the $\Delta^k$ simplex are bounded, the first term tends to zero as $n \to \infty$. We hence find that $
    \CP\left(n \mapsto \frac{1}{n} \sum_{i=1}^n m_i\right) = \CP\left(n \mapsto \frac{1}{n} \sum_{i=1}^n \tilde{m}_i\right)    $.
\end{proof}

\subsubsection{Proof of Proposition~\ref{prop:forwarddir})}
We now prove Proposition~\ref{prop:forwarddir}, for which we require some helpful notation and helper lemmas. Our proof takes some inspiration from \citep[Appendix B]{frohlich2024strictly}.

We define some notation first. Recall that for a set $A$ and some $n \in \mathbb{N}$, we have defined $A^n \coloneqq \{a=(a_1,..,a_n) : a_i \in A\}$. If $a=(a_1,..,a_u) \in A^u$ and $b=(b_1,..,b_v) \in A^v$, we write $a \ccm b \coloneqq (a_1,..,a_u,b_1,..,b_v)$. For $a \in A^u$ and $k\in \mathbb{N}$, we write $a^k \coloneqq a \ccm a \ccm .. \ccm a \in A^{k*u}$. 
Our notation uses different letters for points in the simplex $p,q$ (like relative frequencies), infinite sequences of probability measures $a^\infty,b^\infty$ and finite sequences of probability measures $a,b$ (of which the average is in the simplex). 
For $a=(a_1,..,a_u) \in (\Delta^k)^u$, we denote its average as $r(a) \coloneqq \frac{1}{n} \sum_{i=1}^u a_i$, $r(a) \in \Delta^k$. This average corresponds to the expectation of an independent product of the $a_i$. 
We denote the Euclidean norm by $\|\cdot\|$. As a useful fact, any $p \in \Delta^k$ has $\|p\| \leq 1$. Thus $\|p-q\| \leq 2$ if $p,q$ are in the simplex.
For a point $p\in \Delta^k$ and a line segment $[q,r] \subset \Delta^k$ we say that $p$ is at most $\delta$ away from $[q,r]$ iff $\exists x \in [q,r] : d(p,x)\leq \delta$.

\begin{lemma}[{\citep[Lemma B.9]{frohlich2024strictly}}]
\label{lemma:relfreqsconvex}
    Let $a=(a_1,..,a_u) \in (\Delta^k)^u$, $b=(b_1,..,b_v) \in (\Delta^k)^v$. Then 
    \[
    r(a \ccm b) = \frac{u}{u+v} r(a) + \frac{v}{u+v} r(b) = (1-\alpha) r(a) + \alpha r(b),
    \] where $\alpha \coloneqq v/(u+v)$ and $1-\alpha=u/(u+v)$.
\end{lemma}
\begin{proof} Just compute
    \begin{align}
        r(a \ccm b) &= \frac{1}{u+v} \left(\sum_{i=1}^{u} a_i + \sum_{i=1}^v b_i\right)\\
        &= \frac{u}{u+v} \frac{1}{u} \sum_{i=1}^{u} a_i + \frac{v}{u+v} \frac{1}{v} \sum_{i=1}^v b_i\\
        &= \frac{u}{u+v} r(a) + \frac{v}{u+v} r(b).
    \end{align}
\end{proof}
The next lemma will formalize the intuition: if $a \in (\Delta^k)^u$ is a ``long'' finite sequence of measures and $b \in (\Delta^k)^v$ is relatively ``short'', then the relative frequencies $r(a \ccm b)$ will be close to $r(a)$.
\begin{lemma}[{similar to \citep[Lemma B.17]{frohlich2024strictly}}]
\label{lemma:nottoofar}
    Let $a \in (\Delta^k)^u$ and $b=(b_1,..,b_v) \in (\Delta^k)^v$. Then
    \[
    d(r(a), r(a \ccm (b_1, .. , b_i))) \leq \frac{2v}{u+v}, \quad 1 \leq i \leq v.
    \]
\end{lemma}
\begin{proof}
    Consider first the extreme case $i=v$. We use the previous Lemma~\ref{lemma:relfreqsconvex}. Then
    \begin{align}
        d(r(a),r(a \ccm b)) &= \|r(a)-r(a \ccm b)\|\\
        &=\|(1-u/(u+v) r(a) - v/(u+v) r(b)\|\\
        &=v/(u+v) \|r(a)-r(b)\|\\ 
        &\leq 2v/(u+v),
    \end{align}
    since $\|r(a)-r(b)\|\leq 2$ in the simplex, as any $\|p\|\leq 1$ for $p \in \Delta^k$. For $i<v$ we have that $2i/(u+i) \leq 2v(u+v)$ and thus the bound also holds.
\end{proof}

We will later want to approximate measures $q \in\co(\mM)$ by finite sequences of measures from $\mM$. That this is possible is formalized in the following; the approximation quality of course depends on the length of the finite sequence.
\begin{lemma}[Approximation lemma]
\label{lemma:approximationlemma}
    Let $\emptyset \neq \mathcal{M} \subseteq \Delta^k$ and $q \in\co(\mM)$. Let $v \in \mathbb{N}$ be given. Then there exists some $b \in \mathcal{M}^v$  such that $d(q,r_b) \leq 4(k+1)/v$.
\end{lemma}
\begin{proof}
First note that due to Carathéodory's theorem, we need at most $k+1$ points from $\mathcal{M}$ to express $q \in\co(\mM)$, that is, we can write $q = \sum_{i=1}^w \lambda_i m_i$ where $w \leq k+1$ and $m_i \in \mathcal{M}$, $\sum_{i=1}^w \lambda_i = 1$.
    Write
    \[
    p \coloneqq \sum_{i=1}^w \frac{\lfloor v \lambda_i \rfloor}{v} m_i + \xi p_0, \quad \xi \coloneqq 1 - \sum_{i=1}^w \frac{\lfloor v \lambda_i \rfloor}{v},
    \]
    for some $p_0 \in \mM$. Note that $p \in \Delta^k$ and it is expressed as a convex combination of measures from $\mM$, where all coefficients are rational with common denominator $v$. Implicitly here we constructed a sequence of measures \[
    b \coloneqq \underbrace{m_1 \ccm .. m_1}_{\lfloor v \lambda_1 \rfloor \text{ many times}} \ccm \; .. \ccm \underbrace{m_w \ccm .. m_w}_{\lfloor v \lambda_w \rfloor \text{ many times}} \ccm \underbrace{p_0}_{\xi v \text{ many times}}.
    \]
    Note that $b$ is a finite sequence of measures from $\mM$ of length $v$ and $r_b \in \Delta^k$. It holds that $p=r_b$.

    First, note that if $v \in \mathbb{N}$, then $0 \leq \lambda_i - \lfloor v \lambda_i \rfloor / v \leq 1/v$ since $0 \leq v\lambda_i - \lfloor v\lambda_i \rfloor \leq 1$.
    We then find that 
    \begin{align}
        d\left(\sum_{i=1}^w \frac{\lfloor v \lambda_i \rfloor}{v} m_i, \sum_{i=1}^w \lambda_i m_i\right) &= \left\| \sum_{i=1}^w \left(\lambda_i - \frac{\lfloor v \lambda_i \rfloor}{v} \right) m_i\right\|\\
        &\leq \sum_{i=1}^w \left(\lambda_i - \frac{\lfloor v \lambda_i \rfloor}{v}\right) \|m_i\| \\
        &\leq \frac{2w}{v}.
    \end{align}
    Furthermore, $\xi = 1 - \sum_{i=1}^w \frac{\lfloor v \lambda_i \rfloor}{v} \leq w/v$, which holds because
    we have seen that $\lambda_i - \frac{\lfloor v \lambda_i \rfloor}{v} \leq 1/v$ and $\sum_{i=1}^w \lambda_i = 1$. Therefore
    \begin{align}
        d\left(\sum_{i=1}^w \frac{\lfloor v \lambda_i \rfloor}{v} m_i, \sum_{i=1}^w \frac{\lfloor v \lambda_i \rfloor}{v} m_i + \xi p_0\right) \leq \xi||p_0|| \leq 2\xi \leq \frac{2w}{v}.
    \end{align}
    Taking this together we obtain by the triangle inequality
    \[
    d\left(q, \sum_{i=1}^w \frac{\lfloor v \lambda_i \rfloor}{v} m_i + \xi p_0\right) \leq 4(k+1)/v.
    \]
\end{proof}
The next lemma formalizes the following: if we have already accumulated a finite sequence of measures $a \in (\Delta^k)^u)$, so that our current average $r(a)$ is at most $\delta$ far away from some measure $p$ for sufficiently large $\delta$, then we can find some new measures $(d_1,..,d_{v})$ from $\mM$ so that the new average $r(a \ccm (d_1,..,d_{v})$ is again at most $\delta$ far away from $p$. The smaller $\delta$ is, the more stringent the condition, so $v$ needs to be smaller. The purpose of this lemma is to ``approximately stay where we are'' but grow the length of the finite sequence of measures. This in turn then allows for higher quality approximations as in the previous lemma.
\begin{lemma}[The ``stay and grow'' lemma]
\label{lemma:stayandgrow}
    Let $p \in\co(\mM)$, $q \in B_{\delta}(p)$, $q = r_a, a \in (\Delta^k)^u$, and assume $\delta>0$ is so that it satisfies \[
    \exists v \in \mathbb{N}: \frac{2v}{u+v} \leq \delta \text{ and  } 4(k+1)/v \leq \delta.
    \]
    Then $\exists b_1,..,b_{v} \in \mM $ so that $ r(a \ccm (b_1, .. b_{v})) \in B_\delta(p)$, as well as
    \[
    d(r(a \ccm (b_1, .., b_i)), p) < 2\delta \;\; \forall 0 \leq i \leq v.
    \]
\end{lemma}
\begin{proof}
    Take the above $v \in \mathbb{N}$ and use the approximation lemma~\ref{lemma:approximationlemma} to obtain some $b \in \mM^{v}$-approximation of $p$, call $\tilde{p} \coloneqq r(b)$, so that $d(p,\tilde{p}) \leq \delta$. Then the half-open segment $[q,\tilde{p})$ lies fully in $B_\delta(p)$, since $q$ itself is in $B_\delta(p)$. We can then construct $a \ccm (b_1, .., b_{v})$ and it is guaranteed that $r(a \ccm (b_1, .., b_{v})) \in [q,\tilde{p}) \subset B_\delta(p)$. And since $2v/(u+v)\leq \delta$ by assumption, we have
    \[
    d(r(a \ccm (b_1, .., b_i)), p) \leq d(r(a \ccm (b_1, .., b_i)), q) + d(q,p) < 2v/(u+v) + \delta \leq 2 \delta
    \]
    due to Lemma~\ref{lemma:nottoofar}.
\end{proof}
\begin{lemma}
    Let $\mathcal{N} \subseteq \Delta^k$. Then there exists an $\varepsilon$-cover of $\mathcal{N}$ using finitely many points in $\mathcal{N}$. Formally, $\forall \varepsilon>0 : \exists n_1, .., n_o \in \mathcal{N} : \mathcal{N}_\varepsilon \coloneqq \{n_1, .., n_o\} : \forall n \in \mathcal{N} : \exists n_i \in \mathcal{N}_\varepsilon : d(n,n_i)<\varepsilon$. We call the $n_i$ the $\varepsilon$-centers.
\end{lemma}
\begin{proof}
    Any subset of $\Delta^k$ is precompact, meaning its closure is compact. Then see \eg Remark 4.2.3 in \citep{vershynin2018high} for the result.
\end{proof}
Finally, we need to exploit the connectedness of $\mathcal{N}$.
\begin{lemma}
    Let $\mathcal{N} \subseteq \Delta^k$ be connected. Then $\forall \varepsilon>0$ and any $a,b \in \mathcal{N}$ there exists an $\varepsilon$-chain between $a$ and $b$, that is, a finite set of points $n_1,..,n_o \in \mathcal{N}$ such that $n_1 = a$, $n_o=b$ and $d(n_i,n_{i+1})<\varepsilon$.
\end{lemma}
\begin{proof}
    A connected subset of $\Delta^k$ is a connected metric space by inheriting the Euclidean metric. Any connected metric space is $\varepsilon$-chainable, see \eg \citep[p. 148]{joshi1983introduction}.
\end{proof}
From this is the next lemma readily follows.
\begin{lemma}
\label{lemma:ourepsiloncover}
  Let $\mathcal{N} \subseteq \cobar(\mM) \subseteq \Delta^k$ be connected. Then $\forall \varepsilon>0$ there exists a finite $\varepsilon$-cover of $\mathcal{N}$, which also forms an $\varepsilon$-chain, where all $\varepsilon$-centers are in $\co(\mM)$ (possibly outside of $\mathcal{N}$), but at most $\zeta$-far away from $\mathcal{N}$. Formally, $\forall \varepsilon>0 : \exists p_1,..,p_o \in\co(\mM) : \forall n \in \mathcal{N} : \exists p_i : d(n,p_i)<\varepsilon$; also, $d(p_i,p_{i+1})<\varepsilon$; and $\forall p_i: \exists n \in \mathcal{N}: d(p_i,n)<\zeta$.
\end{lemma}
\begin{proof}
Let $\varepsilon>0$ be arbitrary. Take some $\varepsilon'+2\zeta<\varepsilon$ and use the previous two Lemmata to obtain points $n_1,..,n_0 \in \mathcal{N}$ which form an $\varepsilon'$-cover of $\mathcal{N}$ and also an $\varepsilon'$-chain. Now to each $n_i$ associate some $p_i \in \co \mM$ so that $d(n_i,p_i)<\zeta$, which is possible by assumption that $\mathcal{N} \subseteq \cobar \mM$. By the triangle inequality, these $p_i$ satisfy all desiderata.
\end{proof}
We are now ready to prove Proposition~\ref{prop:mainproptheoretical}.
\begin{proof}
     Recall that we work in the simplex $\Delta^k$. We begin by setting up some infinite sequences of parameterization data, which we then use for a constructive method.
    Let $\varepsilon^\infty=(2,\varepsilon_2,..)$ and $\delta^\infty=(4(k+1),\delta_2,..)$ and $\zeta^\infty=(\zeta_1,..)$ be decreasing and strictly positive sequences satisfying $\varepsilon^\infty \downarrow 0$, $\delta^\infty \downarrow 0$ and $\zeta^\infty \downarrow 0$, to be interpreted as tolerance parameters.
    For each $\delta_i>0$, we compute 
    \[
    v_i \coloneqq \lceil 4(k+1)/\delta \rceil + 1 \in \mathbb{N},
    \]
    which guarantees that $4(k+1)/v_i < \delta_i$. Intuitively, this will be the length with which we invoke the approximation lemma at iteration $i$. However, to move ``safely'' we need to ensure the finite sequence is long enough already. For this, we compute a minimally required length $l_i \in \mathbb{N}$ as the smallest natural number that satisfies $0<2 v_i / (l_i + v_i)\leq \delta_i$. Clearly, by making $l_i$ large enough this can always be guaranteed.
    For example, for $\delta_1=4(k+1)$ we find that $v_1=2$ and $l_i=1$. In addition, for each $i \in \mathbb{N}$, from invoking Lemma~\ref{lemma:ourepsiloncover} with $\varepsilon_i$ and $\zeta_i$, we obtain a finite set of points $E_i = (c_1^i,..,c_{o(i)}^i) \in \co(\mM)^{o(i)}$ which form an $\varepsilon_i$-cover of $\mathcal{N}$ and are also an $\varepsilon_i$-chain, with some finite length $o(i) \in \mathbb{N}$ depending on $i$. For each of these points $c_j \in E_i$, we have $\exists n \in \mathcal{N}: d(c_j,n)<\zeta_i$.

    We now describe how a sequence of measures can be constructed that yields the desired $\mathcal{N}$ as the set of cluster points. The procedure is given as a countable infinity of subsequent iterations, each parameterized by $(\varepsilon_i,\delta_i,\zeta_i,v_i,l_i,E_i)$. 

    \textbf{Iteration $i=1$.}
    For the first iteration $i=1$, we begin by initializing a finite sequence of measures of length $1$ as $m \coloneqq (m_0)$ for some arbitrary measure $m_0$ from $\mM$.

    We now describe how some iteration $i \in \{2,..,\}$, of our algorithm works. We assume that we have already accumulated a finite sequence of measures $m=(m_1,..,m_n)$ from $\mM$. Intuitively, we now want to ``get very close to'' all measures in $E_i$.

\textbf{Iteration $i$}.
At the beginning of iteration $i$ we make two key inductive assumptions:
\begin{enumerate}
    \item The current average $r(m)$ lies in $B_{\delta_i}(c_1^i)$.
    \item  The length $|m|=n$ satisfies $n \geq l_i$.
\end{enumerate}
From the $\varepsilon_i$-centers together with the first point of the $\varepsilon_{i+1}$-centers, create a list of pairs of the form
\[
C = (c_1^i,c_2^i), .., (c_{o(i)-1}^i, c_{o(i)}^i), (c_{o(i)}^i, x_1), .. , (x_t, c_1^{i+1}),
\]
where we also connected $c_{o(i)}^i$ and $c_1^{i+1}$ with an $\varepsilon_i$-chain $c_{o(i)}^i,x_1,..,x_t,c_1^{i+1}$, where all points satisfy that their distance to $\mathcal{N}$ is at most $\zeta_i$ (possible by the reasoning of Lemma~\ref{lemma:ourepsiloncover}). Note that all such points are in $\co(\mM)$, but possibly outside of $\mathcal{N}$. 

\textit{The Pair Subroutine}. 
We now consider each pair $(c_1,c_2) \in C$ with the above order. For each pair, at the beginning of our construction, the current average satisfies $r(m) \in B_{\delta_i}(c_1)$, and we want to append finitely many measures $\tilde{m}=(..)$ from $\mM$ to $m$ so that the new average $r(m\ccm \tilde{m})$ is in $B_{\delta_i}(c_2)$ at the end of the pair subroutine.
To now get closer to $c_2$ we first approximate $c_2$ by $\tilde{c}_2=r_b$ as in the approximation lemma~\ref{lemma:approximationlemma}, using the length $v_i$ for some $b \in (\mM)^{v_i}$, and we get that $d(c_2,\tilde{c}_2) \leq 4(k+1)/v_i < \delta_i$.
Noting that there must exist some $k \in \mathbb{N}$ such that
\begin{equation}
\label{eq:first_implementation_gap}
r(m \ccm b^k) \in B_{\kappa}(\tilde{c}_2),
\end{equation}
where $\kappa>0$ is chosen such that $\kappa + d(c_2,\tilde{c}_2) < \delta_i$, which is possible. The existence of such $k \in \mathbb{N}$ for arbitrarily small $\kappa$ follows from Lemma~\ref{lemma:relfreqsconvex}.
Thus we can append $b$ often enough to end up in the $\delta_i$-ball around $\tilde{c}_2$, since $d(r(m \ccm  b^k),c_2) \leq d(r(m \ccm  b^k), \tilde{c}_2) + d(\tilde{c}_2,c_2) \leq \delta_i$. Now set $\tilde{m} \coloneqq m \ccm  b^k$.

We now observe that the maximum distance we might have been away from $\mathcal{N}$ during this procedure is bounded by $\varepsilon_i+2\delta_i+\zeta_i$: first, we have moved on the $[r(m),\tilde{c}_2]$ line, and potentially went was far as $\delta_i$ away from it during appending but no further. Why?:
\[
d(m \ccm b^{k'}, m \ccm b^{k'} \ccm (b_1,..,b_j)) \leq 2v_i / (l_i + v_i) \leq \delta, \; \; \forall k' \in \{0,..,k-1\} \; \; \forall j \in \{1,..,v_i\},
\]
noting that $2v_i/(\tilde{l}_i+v_i) < 2v_i/(l_i+v_i)$ for any $\tilde{l}_i>l_i$.

 Second, each point $e$ on $[r(m),\tilde{c}_2]$ is at most $\delta_i$ away from $[c_1,c_2]$, \ie there exists a point $c \in [c_1,c_2]$ such that $d(e,c)<\delta_i$. Why? Write $e = (1-\alpha)r(m) + \alpha \tilde{c}_2$, $\alpha \in [0,1]$. Consider the point $f=(1-\alpha)c_1 + \alpha c_2$. Then 
\begin{align}
    d(e,f) &= \| (1-\alpha)r(m) + \alpha \tilde{c}_2 - (1-\alpha)c_1 + \alpha c_2)\\
    &\leq (1-\alpha)\|r(m)-c_1\| + \alpha\|\tilde{c}_2 - c_2\|\\
    &\leq (1-\alpha)\delta_i + \alpha \delta_i\\ 
    &= \delta_i,
\end{align}
since $d(r(m),c_1)<\delta_i$ and $d(\tilde{c}_2,c_2)<\delta_i$. Now, since each point on $[c_1,c_2]$ is at most $\varepsilon_i+\zeta_i$ away from $\mathcal{N}$ due to the $\varepsilon_i$-chain property, we conclude by the triangle inequality that the new averages in this iteration $i$ were never further than $\varepsilon_i+2\delta_i+\zeta_i$ away from $\mathcal{N}$. Now set $m \coloneqq \tilde{m}$ and proceed with the next pair.

\textit{[End of Pair Subroutine]}

After executing the pair subroutine for all pairs, we have accumulated a finite sequence of measures $m$ such that $r(m) \in B_{\delta_i}(c_1^{i+1})$.
Now we would like to move from a $\delta_i$-ball around it to some $\delta_{i+1}<\delta_i$ ball around it, to satisfy the inductive assumption of the next iteration $i+1$.
For this, we invoke the ```stay and grow'' lemma~\ref{lemma:stayandgrow} with $v=v_i$ to end up again within the $\delta_i$-ball around $c_1^{i+1}$ and grow the length of the finite sequence until it has reached at least $l_{i+1}$: applying the lemma once yields some $b \in \mM^{v_i}$ so that 
\[
r(m \ccm b) \in B_{\delta_i}(c_1^{i+1}) \;\; \text{ and } \;\; d(r(m \ccm (b_1,..,b_j)),c_1^{i+1})<2\delta_i, \quad \forall j \in \{1,..,v_i\}.
\]
In words, the lemma guarantees that when making these appendings, the average is never further than $2 \delta_i$ away from $c_1^{i+1}$. Invoking the lemma $r$ many times, for some finite $r \in \mathbb{N}$, set $m \coloneqq m \ccm b^{r}$, guarantees that we reach length $|m|=l_{i+1}$.

We now show that this allows us to safely move to a $\delta_{i+1}$-ball around $c_1^{i+1}$. Our current location in the simplex is $p=r(m)$.
We invoke the approximation lemma with $v=v_{i+1}$ around $c_1^{i+1}$ to obtain a $q=r_b$, $b \in \mM^{v_i+1}$, with $d(c_1^{i+1},q)<\delta_{i+1}$. But then there exists some $\kappa>0$ so that $\kappa + d(c_1^{i+1},q) < \delta_{i+1}$. Now we can move on the line $[p,q]$ until we are in a $\kappa$-ball around $q$, which in turn means we are in a $\delta_{i+1}$-ball around $c_1^{i+1}$. Formally:
\begin{equation}
\label{eq:second_implementation_gap}
\exists k \in \mathbb{N}:  r(m \ccm b^k) \in B_\kappa(q).
\end{equation}
How ``safe'' are these iterative appendings? Since $c_1^{i+1}$ is at most $\zeta_i$-far away from $\mathcal{N}$, $r(m) \in B_{\delta}(c_1^{i+1})$ and for appending $b$ we have the $\delta_{i+1}$-bound,
\[
d(m \ccm b^{k'}, m \ccm b^{k'} \ccm (b_1,..,b_j)) < \delta_{i+1}, \; \; \forall k' \in \{0,..,k-1\} \; \; \forall j \in \{1,..,v_i\},
\]
we have been at most $\delta_i + \delta_{i+1}+\zeta_i<2\delta_i + \zeta_i$ away from $\mathcal{N}$. Now set $m \coloneqq m \ccm b^k$.

To conclude the logic of an iteration $i$, observe that all new averages during this iteration have been at most $\varepsilon_i+2\delta_i+\zeta_i$ away from $\mathcal{N}$, and also we have been $\varepsilon_i + \delta_i + \zeta_i$-close to any point on $\mathcal{N}$.

It remains to argue that the set of cluster points coincides exactly with $\mathcal{N}$. First, obviously the above procedure guarantees that any $n \in \mathcal{N}$ is a cluster point: as $\varepsilon^\infty \downarrow 0$, $\delta^\infty \downarrow 0$ and $\zeta^\infty \downarrow 0$, we enter any neighborhood around any point on $\mathcal{N}$ infinitely often. To see that no other point (outside of $\mathcal{N}$) can be a cluster point, just observe that the averages generated by the above construction get arbitrarily close to $\mathcal{N}$, and $\mathcal{N}$ is already closed by assumption. So for any point $p \notin \mathcal{N}$ we can find an open neighborhood $U$ around it so that $U \cap \mathcal{N} = \emptyset$, which is sufficiently small so that our algorithm never enters this neighborhood, hence the point cannot be a cluster point.
    
\end{proof}

\subsubsection{Concluding the proof of Proposition~\ref{theorem:maintheoremffivanenko}}
It remains only to show the backward direction for Proposition~\ref{prop:mainproptheoretical}.
\begin{proposition}
    Let $\emptyset \neq \mathcal{M} \subseteq \Delta^k$ and $m^\infty=(m_1,m_2,..)$ a sequence of measures from $\mM$. Then the set $\CP\left(n \mapsto \frac{1}{n} \sum_{i=1}^n m_i\right) $ must be closed, connected and a subset of $\cobar(\mM)$.
\end{proposition}
\begin{proof}
    It is clear in general that a set of cluster points must be closed in the respective topology. For connectedness of $\CP\left(n \mapsto \frac{1}{n} \sum_{i=1}^n m_i\right) $, see \citep[Lemma B.10]{frohlich2024strictly}, stated for the case where $\mathcal{M}$ is the set of vertices of $\Delta^k$, but the argument applies for any $\mathcal{M} \subseteq \Delta^k$.
    Finally, assume by contradiction that some point $p \notin \cobar(\mM)$ is a cluster point of $n \mapsto \frac{1}{n} \sum_{i=1}^n m_i$. But for any $n \in \mathbb{N}$, $\left(\frac{1}{n} \sum_{i=1}^n m_i\right) \in \co(\mM)$, and since $p \notin \cobar(\mM)$ we can find an open neighborhood $U$ around it so that $U \cap \cobar(\mM) = \emptyset$. But then the sequence cannot ever enter $U$, and hence $p$ cannot be a cluster point.
\end{proof}

\end{document}